\documentclass[preprint,9pt]{elsarticle}

\usepackage[abs]{overpic}
\usepackage{graphicx}
\usepackage{subfigure}
\usepackage{stmaryrd}
\usepackage{amsfonts, color, amssymb}
\usepackage{amsmath}
\usepackage{epsfig}
\usepackage{psfrag}
\usepackage{paralist}
\usepackage{algorithm}
\usepackage{algorithmicx}
\usepackage{algpseudocode}

\usepackage[misc]{ifsym}

\usepackage[a4paper, total={6in,8in}]{geometry}
\newcommand{\vect}[1]{{\bf #1}}

\newtheorem{remark}{Remark}[section]

\newtheorem{proposition}{Proposition}[section]

\allowdisplaybreaks 

\journal{Journal}

\begin{document}

\begin{frontmatter}

\title{Kernel Free Boundary Integral Method for 3D Stokes and Navier Equations on Irregular Domains}

\author[mymainaddress]{Zhongshu Zhao}
\author[mysecondaryaddress]{Haixia Dong \corref{mycorrespondingauthor}}
\cortext[mycorrespondingauthor]{Corresponding author}
\ead{hxdong@hunnu.edu.cn}
\author[mythirdaddress]{ Wenjun Ying}
\address[mymainaddress]{School of Mathematical Sciences and Institute of Natural Sciences, Shanghai Jiao Tong University, Minhang, Shanghai, 200240, P. R. China.}
\address[mysecondaryaddress]{Corresponding author. MOE-LCSM, School of Mathematics and Statistics, Hunan Normal University, Changsha, Hunan 410081, P. R. China.}
\address[mythirdaddress]{School of Mathematical Sciences, MOE-LSC and Institute of Natural Sciences, Shanghai Jiao Tong University, Minhang, Shanghai 200240, P. R. China.} 

\begin{abstract}
 A second-order accurate kernel-free boundary integral method is presented for Stokes and Navier boundary value problems on three-dimensional irregular domains. It solves equations in the framework of boundary integral equations, whose corresponding discrete forms are well-conditioned and solved by the GMRES method. A notable feature of this approach is that the boundary or volume integrals encountered in BIEs are indirectly evaluated by a Cartesian grid-based method, which includes discretizing corresponding simple interface problems with a MAC scheme, correcting discrete linear systems to reduce large local truncation errors near the interface, solving the modified system by a CG method together with an FFT-based Poisson solver. No extra work or special quadratures are required to deal with singular or hyper-singular boundary integrals and the dependence on the analytical expressions of Green's functions for the integral kernels is completely eliminated. Numerical results are given to demonstrate the efficiency and accuracy of the Cartesian grid-based method.

\end{abstract}

%


\end{frontmatter}


\section{Introduction}
The Stokes and Navier problems are two important models constructed in incompressible fluid and solid mechanics, and have wide applications in engineering and sciences, such as lubrication theory \cite{silva2016stokes}, porous media flow \cite{layton2002coupling},  tissue engineering \cite{wei2010differential}, biomedical science \cite{wei2013multiscale} and so on. Therefore, it is always of great interest to find simple, effective and robust numerical schemes for solving these models.

For such problems defined on irregular or complex geometries, a traditional numerical method such as the finite element method (FEM) with a body-fitted grid suffers from a major challenge on efficient mesh generation and accurate solution of corresponding systems. Especially in three dimensions or when moving boundaries are involved, mesh generation and re-meshing become difficult and time-consuming. In addition, another difficulty is the design of robust and fast solvers for the resulting discrete equations. Although multilevel techniques such as multi-grid or domain decomposition have been extended to unstructured grids, vector PDEs or indefinite operators like the Stokes equations have not been widely applied.

Also, many numerical methods, such as the finite difference method, do not naturally apply to unstructured meshes. In order to avoid these drawbacks, the use of Cartesian grid-based methods has become quite widespread. Representative numerical methods of this type partially include immersed boundary method (IBM) \cite{peskin1977numerical}, immersed interface method (IIM) \cite{leveque1994immersed}, CutFEM \cite{hansbo2014cut}, extended finite element method (XFEM)\cite{chessa2003extended}, Nitsche's XFEM \cite{wang2019nonconforming}, immersed finite element method (IFEM) \cite{he2011immersed}, matched interface and boundary (MIB) method \cite{zhou2006high} and so on.  To maintain the desired accuracy, techniques such as smoothing or regularization of discontinuities, correction of the discretization schemes and modification of the approximation functions or basis are usually employed. 
Most Cartesian grid-based methods enable one to employ much simple meshes, but in some cases, fast methods are not straightforward to apply.

Boundary integral methods (BIMs) have been used most extensively in the case of ellipses because they have the significant advantages of handling complex or irregular domains and using fast algorithms to dramatically decrease the computational cost.  The main idea is to embed the complex or irregular domain into a larger regular domain, and then the boundary value problems can be reformulated into Fredholm BIEs of the second kind, which leads to the fact that only the domain boundary or surface is needed to be partitioned, thus avoiding the generation of high-quality boundary-fitted mesh on irregular domains, and considerably reducing the dimensionality of unknowns in the solution.  After the first numerical implementation of the boundary integral formulation for Stokes flow reported by Youngren and Acrivos in \cite{youngren1975stokes}, BIMs have played an important role in fluid mechanics, elasticity and other application areas \cite{greengard1996integral,jiang2012integral,rachh2016integral,abousleiman1994boundary,zeb1998boundary,grigoriev2005multi,biros2004fast}. However, there exist several  potential practical issues that have prevented the broad application of approaches of this type. For example, the singularities of the fundamental solution can involve increased computational costs and implementation complexity when computing the near field interaction. Although various fast techniques have been developed to speed up the calculation of singular integrals (see e.g.\cite{huang1993some,beale2001method,klaseboer2012non,tlupova2013nearly} and references therein for further information), it is still an active research topic. In addition, the unavailability of the analytical expression of the kernel functions also restricts the traditional BIMs to the constant coefficients boundary value problems in the free space subject to the far field radiation condition or in a rectangle domain subject to the periodic boundary condition.  

For these reasons, a different approach has emerged. This method, referred to as Kernel-free boundary integral (KFBI) method, is a generalization of the traditional BIMs, particularly the grid-based BIM by Mayo et al. \cite{mayo1984fast, mayo1985fast,mayo1992rapid} and J.T.Beale et al \cite{beale2004grid}. It evaluates boundary and volume integrals indirectly by a Cartesian grid-based method, thus possessing the following two most prominent features: i) it does not require the explicit expressions of Green's function or special quadratures formulas to directly calculate integrals, especially nearly singular or hyper-singular boundary integrals, so that the dependence on the kernel can be completely eliminated in practice; ii) it reformulates the boundary value problems as the Fredholm BIE of the second kind, helping to eliminate the ill-conditioning property of the original problems so that the number of Krylov subspace iterations is essentially independent of the discretization parameter or system dimension. The KFBI method has been developed to be a general method for two-dimensional elliptic PDEs \cite{ying2007kernel,ying2013kernel,ying2014kernel,ying2018cartesian,xie2019fourth,xie2020high,cao2022kernel}, but for three-dimensional problems, it is now under intensive development.

This work extends the KFBI method to solve Stokes and Navier equations in three-dimensional irregular domains. Based on potential theory, the solution of a Dirichlet problem is written as the sum of a volume integral and a double-layer boundary integral with an unknown density. These integrals are evaluated indirectly by a Cartesian grid-based method, which primarily consists of two steps: (1) solving corresponding equivalent but simple interface problems in an extended cubic region, (2) extracting the boundary value of the integrals by a procedure of polynomial interpolation. During the calculation, the equivalent Stokes and Navier interface problems are discretized in a uniform mixed formulation with a modified MAC Scheme, generalized slightly by allowing a pressure term in the continuity equation.  The resulting linear system is solved efficiently by the CG method together with an FFT-based Poisson solver. The Cartesian grid-based indirect evaluation technique has the superiorities of requiring no extra work or special quadratures to handle singular or hyper-singular boundary integrals without the need to analytical expressions of Green's functions for the integral kernels.

In addition, no unstructured triangulation of the surface is required in the KFBI method. It only uses some quasi-uniform control points, which are represented by intersection points of the surface with an underlying Cartesian grid, to discretize the density and the boundary integral equations. Such a selection of control points makes the interpolation stencils in the integral evaluation convenient to choose and locally uniform on a coordinate plane in three space dimensions. As the intersection points of an implicit surface with Cartesian grid line can be found straightforwardly in three space dimensions, it is very easy to implement the algorithm. These attributes have special importance for the time-dependent problems with moving boundaries. Numerical results show that the KFBI method is efficient and accurate in handling incompressible fluid and solid mechanics problems on irregular domains.

The remainder of this article is organized as follows. In section 2, the Stokes and Navier problems are described respectively.  The main idea of the KFBI algorithm is given in section 3. 
The essential implementation details for integral evaluation are given in section 4. In section 5,  numerical results are provided to validate the effectiveness of the proposed method. Concluding remarks and some discussions are put in section 6. 

\section{Boundary Value Problems}
Let $\Omega\in \mathbb{R}^3$ be a  bounded domain with smooth boundary $\partial \Omega$, which is in general irregular and complex. The steady-state incompressible Stokes equations  considered in this work are given by 
\begin{subequations}
\label{SP}
\begin{align}
	 -\mu_s \Delta\vect   u_s+\nabla p_s&= \vect   f_s, \quad\, \hbox{in}\; \Omega,\label{SP-1}\\
\nabla\cdot \vect   u_s&= 0,\quad\;\, \hbox{in}\; \Omega,\label{SP-2}\\
\vect   u_s&=\vect   g_s, \quad\hbox{on}\; \partial \Omega,\label{SP-3}
\end{align}
\end{subequations}
where  $\vect u_s = (u^{(1)}_s, u^{(2)}_s, u^{(3)}_s)^T$ stands for the velocity vector, 
$p_s$ represents the pressure,
$\mu_s$ is the fluid viscosity coefficient, $\vect f_s=(f_s^{(1)}, f_s^{(2)}, f_s^{(3)})^T$ denotes an external force and $\vect g_s = (g_s^{(1)}, g_s^{(2)}, g_s^{(3)})^T$ is the Dirichlet-type boundary condition of $\vect u_s$ on the boundary $\partial \Omega$. Assume that $\mu_s$ is a constant function on $\Omega$ and $\vect g_s$ satisfies the compatibility condition 
$$\int_{\partial \Omega} \vect g_s\cdot \vect n ds= 0,$$ 
where $\vect n$ denotes the unit outward normal vector on $\partial\Omega$.
The 3D Navier equations are also considered here, which are governed by
\begin{subequations}
\label{Nm}
\begin{align}
\nabla\cdot\pmb \sigma(\vect u_e)+\vect f_e&=\vect 0, \,\;\;\quad{\rm in}\; \Omega,\label{Nm1}\\
\vect u_e &= \vect g_e, \;\quad{\rm on}\;\partial\Omega,\label{Nm2}
\end{align}
\end{subequations}
where $\vect u_e=(u^{(1)}_e, u^{(2)}_e, u^{(3)}_e)^T$ is the displacement field, $\vect f_e=(f^{(1)}_e, f^{(2)}_e, f^{(3)}_e)^T$ is  a given body force,  $\vect g_e = (g^{(1)}_e, g^{(2)}_e, g^{(3)}_e)^T $ represents  the displacement on the boundary. The stress tensor is given by
\begin{align}
\label{Strain}
\pmb \sigma(\vect u_e)=\lambda\nabla\cdot\vect u_e\vect I+2\mu\pmb \epsilon(\vect u_e),
\end{align}
where $\pmb  \epsilon(\vect u_e)=\dfrac{1}{2}\big(\nabla \vect u_e+(\nabla \vect u_e)^T\big)$ is the linear strain and $\vect I$ is the $3\times 3$ identical matrix. Furthermore, $\lambda$ and $\mu$ represent the Lam{\'e} coefficients, satisfying
  $$\lambda = \dfrac{E\nu}{(1+\nu)(1-2\nu)}, \hbox{\qquad and \qquad}\mu=\dfrac{E}{2(1+\nu)},$$
with $E$ the Young's modulus and $\nu$ the Poisson's ratio.  By inserting the strain tensor \eqref{Strain} into the equilibrium equation \eqref{Nm1}, the 3D Navier problem can be reformulated as a system of partial differential equations where the unknown function is the displacement field $\vect u_e$:
\begin{equation}
\label{Navier}
-\mu\Delta\vect u_e - (\lambda+\mu)\nabla(\nabla\cdot \vect u_e) = \vect f_e, \;\;\hbox{in}\;\Omega.
\end{equation}

\section{The Kernel Free Boundary Integral Method}

This section gives details of the Cartesian grid-based KFBI method for solving Stokes and Navier problems in three-dimensional irregular domains. 

 To this end, the irregular domain $\Omega$ is embedded into a larger cube $\mathcal{B}$, thus the domain boundary becomes an interface, which separates the cuboid into two disconnected subdomains $\Omega$ and $\Omega^c$. Here, $\Omega^c = \mathcal{B}\backslash \bar{\Omega}$ is the complement of $\Omega$ in $\mathcal{B}$. In the remainder of this article, the boundary $\partial\Omega$ is redefined as $\Gamma$. Next, problems \eqref{SP} and \eqref{Nm} will be reformulated into boundary integral equations respectively. Let $(\vect G_{\vect v}(\vect x, \vect y), G_q(\vect x, \vect y))$ be Green's function pair associated with the Stokes equation \eqref{SP}  that satisfies
\begin{align}
\label{Green-S}
\begin{split}
-\mu_s\Delta \vect G_{\vect v}(\vect x, \vect y) + \nabla G_q(\vect x, \vect y) &= \vect I \delta(\vect x - \vect y), \quad \hbox{in}\; \mathcal{B},\\
\nabla \cdot \vect G_{\vect v}(\vect x, \vect y) &= 0, \qquad\qquad\;\;\,\hbox{in}\; \mathcal{B},\\
\vect G_{\vect v}(\vect x, \vect y) &= \vect 0, \qquad\qquad\;\; \hbox{on}\; \partial \mathcal{B},
\end{split}
\end{align}
for each fixed $\vect y\in \mathcal{B}$.  Let $\vect G(\vect x, \vect y)$ be Green's function associated with Navier system \eqref{Nm} that satisfies 
\begin{equation}
\label{Green-E}
\begin{split}
-\mu\Delta\vect G(\vect x,\vect y)-(\lambda+\mu)\nabla\big(\nabla\cdot \vect G(\vect x, \vect y)\big) &= \vect I\delta(\vect x-\vect y), \;\; \;{\rm in}\;\mathcal{B},\\
\vect G(\vect x, \vect y) &=\vect 0, \qquad\,\qquad\;{\rm on}\;\partial\mathcal{B},
\end{split}
\end{equation}
for each fixed $\vect y\in\mathcal{B}$. Here, the matrix $\vect  I$ denotes the unit matrix in $\mathbb{R}^3$ and $\delta(\vect x-\vect y)$ is the 3D Dirac delta function. All differentiations are carried out with respect to the variable $\vect x$. It is important to point out that the Green's function pair $(\vect G_{\vect v}(\vect x, \vect y), G_q(\vect x, \vect y))$ and Green's function $\vect G(\vect x, \vect y)$ defined in the bounded domain $\mathcal{B}$ are different from the fundamental solution in the free space \cite{hsiao2008boundary}. Their expressions are in general not analytically known, but their existence is guaranteed.

In terms of the Green's function pair $(\vect G_{\vect v}, G_q)$, the solutions $\vect u_s(\vect x)$ and $p_s(\vect x)$ to the Stokes problem \eqref{SP} can be expressed as a sum of the double layer potential and volume integral
\begin{align}
\label{s-S}
	&\vect u_s(\vect x) = \mathcal{G}_{\vect v}\vect f_s(\vect x) - \mathcal{M}_{\vect v}\pmb\varphi_s(\vect x) , \quad \vect x \in \Omega, \\
	&\,p_s(\vect x) = \mathcal{G}_q\vect f_s(\vect x) - \mathcal{M}_q\pmb\varphi_s(\vect x) , \;\quad \vect x \in \Omega,
\end{align}
with density $\pmb\varphi_s(\vect x)$ satisfying the boundary integral equation 
\begin{align}
\begin{split}
\label{BIE-S}
	\frac{1}{2}\pmb\varphi_s(\vect x) - \mathcal{M}_{\vect v}\pmb\varphi_s(\vect x)  = \vect g_s(\vect x)  - \mathcal{G}_{\vect v}\vect f_s(\vect x) , \qquad \vect x \in \Gamma.
\end{split}
\end{align}
Here, the double layer boundary integrals $\mathcal{M}_{\vect v}\pmb\varphi_s$, $\mathcal{M}_q\pmb\varphi_s$ and volume integrals $\mathcal{G}_{\vect v}\vect  f_s$, $\mathcal{G}_q\vect f_s$ are given respectively by 
\begin{align*}
\begin{split}
(\mathcal{M}_{\vect  v}\pmb\varphi_s)(\vect  x) &= \int_{\Gamma} T_s(\vect  G_{\vect  v}, G_q) \pmb\varphi_s(\vect  y) ds_{\vect  y}, \quad\;\;
(\mathcal{G}_{\vect  v}\vect  f_s)(\vect  x) = \int_{\Omega}\vect  G_{\vect  v}(\vect  x, \vect  y)\vect  f_s(\vect  y)d\vect  y,\\
(\mathcal{M}_q\pmb\varphi_s)(\vect  x) &= 2\int_{\Gamma} \dfrac{\partial G_q(\vect  x, \vect  y)}{\partial \vect  n_{\vect  y}} \pmb\varphi_s(\vect  y)ds_{\vect  y},\quad
(\mathcal{G}_q\vect  f_s)(\vect  x) = \int_{\Omega} G_q(\vect  x, \vect  y)\cdot\vect  f_s(\vect  y)d\vect  y,
\end{split}
\end{align*}
with the traction $\vect T_s(\vect u_s, p_s) =- p_s \vect n + \mu (\nabla \vect u_s+\nabla \vect u_s^T) \vect n.$
Similarly, in terms of Green's function $\vect G(\vect x, \vect y)$, the solution $\vect u_e$ to the Navier equations \eqref{Nm} can also be expressed as a sum of the double layer potential and volume integral
\begin{equation}
\label{s-E}
\vect u_e(\vect x) = \mathcal{G}\vect f_e(\vect x) -\mathcal{M}\pmb \varphi_e(\vect x), \quad \vect x\in\Omega,
\end{equation} 
with density $\pmb\varphi_e(\vect x)$ satisfying the boundary integral equation 
\begin{equation}
\label{BIE-E}
\frac{1}{2}\pmb \varphi_e(\vect x) - \mathcal{M}\pmb \varphi_e(\vect x) =\vect g_e(\vect x) - \mathcal{G}\vect f_e(\vect x).
\end{equation}
Here, the double layer boundary integral $\mathcal{M}\pmb \varphi_e$ and volume integral $\mathcal{G}\vect f_e$ are defined respectively by
\begin{equation}
\label{operatorEDG}
(\mathcal{M}\pmb \varphi_e)(\vect x) =\int_{\Gamma}\vect T_{\vect y}\pmb G(\vect x,\vect y)\pmb \varphi_e(\vect y)ds_{\vect y},\qquad
(\mathcal{G}\vect f_e)(\vect x) =\int_{\Omega}\vect G(\vect x,\vect y)\vect f_e(\vect y)d\vect y, 
\end{equation}
with the traction defined by
\begin{equation*}
\vect T(\vect u_e)  =\lambda(\nabla\cdot\vect u_e)\vect n+2\mu\dfrac{\partial \vect u_e}{\partial\vect n}+\mu\vect n\times(\nabla\times \vect u_e).
\end{equation*}
It is pointed out that the subscript $\vect y$ in $\vect T_{\vect y}\vect G(\vect x,\vect y)$ denotes the differentiations in \eqref{operatorEDG} with respect to the  variable $\vect y$. 
Note that the BIEs \eqref{BIE-S} and \eqref{BIE-E} are Fredholm integral equations of the second kind \cite{kress1989linear},  which indicate that the iterative methods, such as the generalized minimal residual (GMRES) method\cite{saad1986gmres}, or even the simple Richardson iteration, converge for each $\vect g_s/\vect g_e$ and  initial data $\pmb \varphi_s^0/\pmb \varphi_e^0$ to the unique solution of \eqref{BIE-S} or \eqref{BIE-E}.

Once the iteration converges, one can get approximation of $\vect u_s$ or $\vect u_e$ respectively according to the representation formula \eqref{s-S} and \eqref{s-E}.
As the Green's function pair $(\vect G_{\vect v},G_q)$ and Green's function $\vect G$ are defined in a bounded domain $\mathcal{B}$, their analytical expressions are un-available, thus it is impossible to directly calculate the boundary and volume  integrals encountered in  BIEs \eqref{BIE-S} and \eqref{BIE-E}.

 For this reason, a Cartesian grid-based KFBI method that calculates integrals indirectly is adopted here, so it is not necessary to know the analytical expression of Green's function.  The main idea of KFBI consists three main steps: 1) The irregular domain $\Omega$ is embedded into a larger cuboid area $\mathcal{B}$, which is easy to obtain a uniform Cartesian grid, thus avoids generating unstructured grids for complex domains effectively. 
Then the boundary value problems \eqref{SP}-\eqref{Nm} are reformulated into Fredholm BIEs of the second kind \eqref{BIE-S} and \eqref{BIE-E}. 2) Evaluation of the integrals encountered in the BIEs are made indirectly by a Cartesian grid-based method, including discretizing the corresponding interface problem with a MAC scheme, correcting the established linear system to reduce the large local truncation errors near the interface, solving the modified system by a CG method together with the FFT-based Poisson solvers,  approximating values of the integrals or its normal flux on the interface by quadratic polynomial interpolation. 3) The BIEs are well-conditioned and the corresponding discrete form can be solved efficiently with a Krylov subspace method, such as the GMRES method, with the number of iterations independent of the mesh size.  At last, the algorithm is summarized as follows:
\begin{algorithm}[!htb]
  \caption{ KFBI method for 3D Stokes and Navier Equations.}
  \label{alg:Framwork}
  \begin{algorithmic}
    \renewcommand{\algorithmicrequire}{ \textbf{1:}}
    \Require {\em Find quasi-uniform control points to discretize BIEs.}
    \begin{itemize}
    \item Set up five different uniform staggered grids with size $h$ covering larger cuboid area $\mathcal{B}$. 
    \item Each grid point is marked as an interior or exterior point according to the interface.
    \item Identify regular and irregular grid nodes of the uniform grids.
    \item Find the intersection points of the boundary and the grid lines, and compute tangential and normal unit directions of the boundary curve at those intersected points.
    \end{itemize}
    \renewcommand{\algorithmicrequire}{ \textbf{2:}}
    \Require {\em Evaluate the  boundary or volume integral at the boundary $\Gamma$.}
    \begin{itemize}
    \item Write the equivalent interface problems of Navier and Stokes as a unified form by introducing an auxiliary unknown pressure. 
    \item Discretize the equivalent interface problems with a finite-difference-based MAC scheme.
    \item Compute jumps of the solutions and their partial derivatives at intersection points.
    \item Correct the right-hand side of the MAC scheme at irregular grid nodes.
    \item Solve the modified linear system by the CG method together with an FFT-based Poisson solver.
    \item Extract the boundary data at intersected points by quadratic polynomial interpolation.
    \end{itemize}
    \renewcommand{\algorithmicrequire}{ \textbf{3:}}
    \Require {\em Make the GMRES iteration with the KFBI method. }
    \begin{itemize}
        \item Evaluate the volume integral boundary data with step 2.
        \item Start the GMRES iteration with the trivial zero initial guess and set up a tolerance.
        \item Evaluate the double layer integral boundary data with step 2.
        \item Update the unknown density $\pmb\varphi_s$ by the GMRES iteration until the residual is smaller than the prescribed tolerance in some norm.
    \end{itemize}
  \end{algorithmic}
\end{algorithm}

Noted that the points used to represent the boundary and discretize the BIE are chosen as intersection points of the boundary surface with the grid lines of the underlying staggered grid, which was originally prescribed in \cite{ying2013kernel} by Wenjun Ying and Wei-Cheng Wang. With this technique, the discretization points are convenient to locate. Moreover, the points are classified as the primary and secondary points and each of them is associated with a component of boundary data during the solution of the BIE, which makes it easy to find out compact and locally uniform  interpolation stencils for boundary interpolation. And an additional equilibration process further ensures  the stability and efficiency of the numerical differentiation when calculating the tangential derivatives of the boundary.  One may refer to the reference \cite{ying2013kernel} for a detailed presentation.

The main challenge in the above algorithm is the evaluation of integrals appeared in BIEs, as the conceivable unavailability of the analytical expressions of Green's functions. To this end, each boundary and volume integral is calculated by interpolating structured grid-based solutions, which avoids directly discretizing the integrals by numerical quadratures, so that the dependence on the analytical expressions of Green's functions can be completely eliminated in implementation. Details of the evaluation are presented in the ensuing sections.

\section{Evaluation of Boundary or Volume Integrals}
In this section, a Cartesian grid-based method for indirectly evaluating boundary and volume integrals is presented. In this method, analytical expressions of Green's functions are no longer needed and the integrals are calculated indirectly through the following steps :
\begin{compactitem}
\item[1).] transforming the boundary or volume integral into equivalent interface problems;
\item[2).] solving the equivalent but simple interface problem under Cartesian mesh;
\item[3).]interpolating the discrete solution on the Cartesian mesh to extract values of the integrals at discretization points of the interface.
\end{compactitem}

\subsection{Equivalent Simple Interface Problems}
The equivalent simple interface problems for the volume and double layer boundary integrals  will be illustrated here, respectively. Let $\vect v^{+}(\vect x)$ and $\vect v^{-}(\vect x)$ be the restrictions of $\vect v(\vect x)$ from the subdomain $\Omega$ and $\Omega^c$. For $\vect x\in \Gamma$, $\vect v^+(\vect x)$ and $\vect v^-(\vect s)$ are interpreted as the limit values of $\vect v(\vect x)$ from the corresponding side of the domain boundary. The jump across the interface $\Gamma$ is denoted by 
\begin{align*}
	\begin{split}
		[\![ \vect v(\vect x) ]\!] = \vect v^{+}(\vect x) - \vect v^{-}(\vect x),\quad\vect x\in\Gamma. 
	\end{split}
\end{align*}

\begin{proposition}\rm
\label{FBI}
For a given function $\vect f_s$ defined on $\Omega$,  the volume integrals $\vect u=(\mathcal{G}_{\vect v}\vect f_s)(
\vect x)$ and $p = (\mathcal{G}_q\vect f_s)(\vect x)$ are generalized solution pair to the following simple interface problem
\begin{align}
\label{InterfaceVolume}
\begin{split}
-\mu_s \Delta\vect   u+\nabla p &= \widetilde{\vect f_s} = \left\{
\begin{aligned}
&\vect f_s, \quad\vect x\in \Omega,\\
&\vect 0, \quad\;\vect x\in \Omega^c,
\end{aligned}
\right. \\
\nabla\cdot \vect   u&=0,\,\;\;\;\vect x\in \Omega\cup\Omega^c,\\
[\![ \vect   u ]\!]&=\pmb 0, \,\;\;\;\vect x\in \Gamma,\\
[\![\pmb \sigma_s(\vect   u, p)\vect   n]\!]&=\pmb 0, \;\;\;\,\vect x\in \Gamma,\\
\vect   u&=\vect   0,  \;\;\;\,\vect x\in \partial\mathcal{B}.
\end{split}	
\end{align}
The interface conditions above imply the continuous property of  the volume potential $\vect u$ as well as its stress tensor $\pmb \sigma_s (\vect u, p)\vect n=(-p\vect I + \mu_s(\nabla \vect u + (\nabla \vect u)^T))\vect n$.
\end{proposition}

\begin{proposition}\rm
\label{DBI}
For a given density function $\pmb \varphi_s$ defined on $\Gamma$, the double layer boundary integrals $\vect u(\vect x) = -(\mathcal{M}_{\vect v}\pmb\varphi_s)(\vect x)$ and $q = -(\mathcal{M}_q\pmb \varphi_s)(\vect x)$  are generalized solution pair to the following simple interface problem
\begin{align}
\label{InterfaceDouble}
\begin{split}
-\mu_s \Delta\vect   u+\nabla p&= \vect 0,\,\,\;\;\vect x\in \Omega\cup\Omega^c,\\
\nabla\cdot \vect   u&=0,\,\;\;\; \vect x\in \Omega\cup\Omega^c,\\
[\![ \vect  u ]\!]&=\pmb \varphi_s, \;\,\vect x\in \Gamma,\\
[\![ \pmb \sigma_s(\vect   u, p)\vect   n]\!]&=\pmb 0,\;\; \,\;\vect x\in \Gamma,\\
\vect   u&=\vect   0,  \;\;\,\,\vect x\in \partial\mathcal{B}.
\end{split}	
\end{align}
The stress tensor $\pmb\sigma_s(\vect   u, p)\vect   n$ is continuous across the interface $\Gamma$ and the double layer potential $\vect u$ has a jump $\pmb \varphi_s(\vect x)$, i.e.
\begin{align*}
	&\vect u^{+}=\frac{1}{2}\pmb \varphi_s - \mathcal{M}_{\vect v}\pmb \varphi_s,  \;\;\;\vect x\in \Gamma,\\
	&\vect u^{-}=-\frac{1}{2}\pmb \varphi_s - \mathcal{M}_{\vect v}\pmb \varphi_s,  \;\;\;\vect x\in \Gamma.
\end{align*}
\end{proposition}

\begin{proposition}\rm
\label{FBIN}
For a given function $\vect f_e$ defined on $\Omega$,  the volume integral $\vect u=(\mathcal{G}\vect f_e)(\vect x)$  is a generalized solution to the following simple interface problem
\begin{equation}
\label{IP1N}
\begin{split}
-\mu\Delta\vect u - (\lambda+\mu)\nabla(\nabla\cdot \vect u) &= \widetilde{\vect f}_e:=\begin{cases}  
\vect f_e(\vect x), \;\vect x\in \Omega,\\
\vect 0, \quad\;\;\;\vect x\in \Omega^c,
\end{cases}\\
[\![ \vect  u ]\!] &= \vect 0, \;\;\;\vect x\in \Gamma,\\
[\![ \pmb  \sigma(\vect u)\vect n ]\!] &= \vect 0, \;\;\;\vect  x\in \Gamma,\\
\vect v&= \vect 0, \;\;\;\vect x\in \partial\mathcal{B}. 
\end{split}
\end{equation}
The interface conditions above imply the continuous property of  the volume potential $\vect u$ as well as its traction $\pmb \sigma(\vect u)\vect n$.
\end{proposition}
\begin{proposition}\rm
\label{DBIN}
For a given density function $\pmb \varphi_e$ defined on $\Gamma$, the double layer boundary integral $\vect u=-\mathcal{M}\pmb \varphi_e$
 is a solution to the following simple interface problem
\begin{subequations}
\label{IP2N}
\begin{align}
-\mu\Delta\vect u - (\lambda+\mu)\nabla(\nabla\cdot \vect u) &= \vect 0,\,\;\;\;\vect x\in \Omega\cup\Omega^c,\label{IP2-1}\\
[\![ \vect   u ]\!] &=\pmb \varphi_e,\;\,\vect x\in \Gamma,\label{IP2-2}\\
[\![ \pmb  \sigma(\vect u)\vect n ]\!] &=\vect 0, \,\;\;\;\vect x\in \Gamma,\label{IP2-3}\\
\vect u&= \vect 0,\;\;\;\,\vect x\in \partial\mathcal{B}\label{IP2-4}. 
\end{align}
\end{subequations}
The discontinuity properties of the double layer potential $\vect u$ imply that 
\begin{subequations}
\begin{align}
	&\vect u^{+}=\frac{1}{2}\pmb \varphi_e - \mathcal{M}_{\vect v}\pmb \varphi_e,  \;\;\; \vect x\in \Gamma,\label{IP2-5}\\
	&\vect u^{-}=-\frac{1}{2}\pmb \varphi_e - \mathcal{M}_{\vect v}\pmb \varphi_e,  \;\;\; \vect x\in \Gamma.\label{IP2-6}
\end{align}
\end{subequations}
\end{proposition}

\begin{remark}
To evaluate the volume or boundary integral in \eqref{BIE-S}, one can turn to solve the interface problems  \eqref{InterfaceVolume} or \eqref{InterfaceDouble}, which can be rewritten into a unified form as,
\begin{subequations}
\label{InterfaceUnfied}
\begin{align}
-\mu_s \Delta\vect   u+\nabla p &= \widetilde{\vect f}_s,\;\;\,\vect x\in \Omega\cup\Omega^c,\label{IUS-1}\\
\nabla\cdot \vect   u&=0,\,\;\;\;\vect x\in \Omega\cup\Omega^c,\label{IUS-2}\\
[\![ \vect  u ]\!]&=\pmb \varphi_s, \;\,\vect x\in \Gamma,\label{IUS-3}\\
[\![ \pmb \sigma_s(\vect   u, p)\vect   n]\!]&=\pmb 0,\;\;\;\,\vect x\in \Gamma,\label{IUS-4}\\
\vect   u&=\vect   0, \;\;\;\,\vect x\in \partial\mathcal{B}.\label{IUS-5}	
\end{align}
\end{subequations}
By the linearity of the problems, the solution to the interface problem \eqref{InterfaceUnfied} is the sum of the solutions to the previous two interface problems \eqref{InterfaceVolume} and \eqref{InterfaceDouble}. Similarly, 
To evaluate the volume or boundary integral in \eqref{BIE-E},  one can turn to solve the interface problems  \eqref{IP1N} or \eqref{IP2N}, which can also be  rewritten into a  unified form as 
\begin{equation}
\label{I}
\begin{split}
-\mu\Delta\vect u - (\lambda+\mu)\nabla(\nabla\cdot \vect u) &= \tilde{\vect f}_e,\;\;\vect x\in \Omega\cup\Omega_c,\\
[\![ \vect   u ]\!] &=\pmb \varphi_e, \;\vect x\in \Gamma,\\
[\![ \pmb  \sigma(\vect u)\vect n ]\!] &=\pmb 0, \;\;\;\vect x\in \Gamma,\\
\vect u&= \vect 0,\;\;\;\vect x\in \partial\mathcal{B}. 
\end{split}
\end{equation} 
By introducing an additional unknown $p=-(\lambda+\mu)\nabla \cdot \vect u$,  the above interface problem \eqref{I} can be transcribed in the form of the following system, 
\begin{align}
\label{I-E}
\begin{split}
-\mu\Delta\vect   u+\nabla p &= \tilde{\vect f}_e,\;\;\,\vect x\in \Omega\cup\Omega^c,\\
\nabla\cdot \vect   u+cp&=0,\,\;\;\;\vect x\in \Omega\cup\Omega^c,\\
[\![ \vect  u ]\!]&=\pmb \varphi_e, \;\,\vect x\in \Gamma,\\
[\![ \pmb \sigma(\vect   u, p)\vect   n]\!]&=\pmb 0,\;\;\;\,\vect x\in \Gamma,\\
\vect   u&=\vect   0, \;\;\;\,\vect x\in \partial\mathcal{B},
\end{split}	
\end{align}
with $c=1/(\lambda+\mu)$.  
It is noted that the system \eqref{I-E} corresponds to the Stokes system \eqref{InterfaceUnfied}.  Thus, the following subsection will mainly focus on numerically solving problems \eqref{InterfaceUnfied}.
\end{remark}

\subsection{Modified Marker-and-Cell Scheme} 
\begin{figure}[!ht]
\centering
\subfigure[Staggered grid]{
\includegraphics[width=0.25\textwidth]{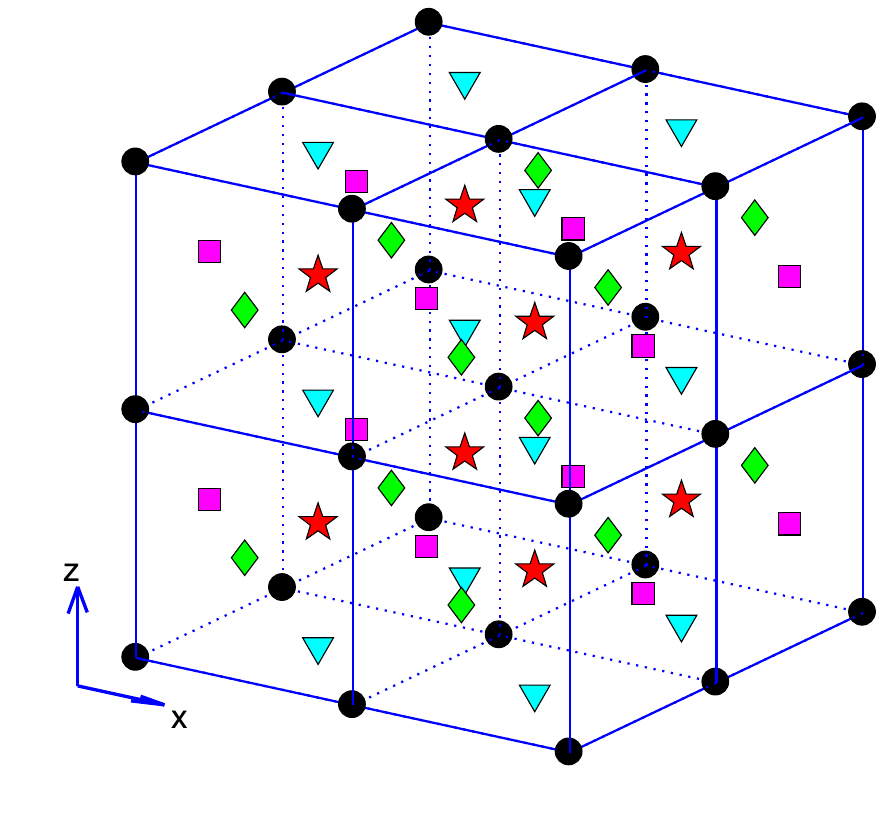}
}
\subfigure[$\mathcal{T}_h$]{
\includegraphics[width=0.25\textwidth]{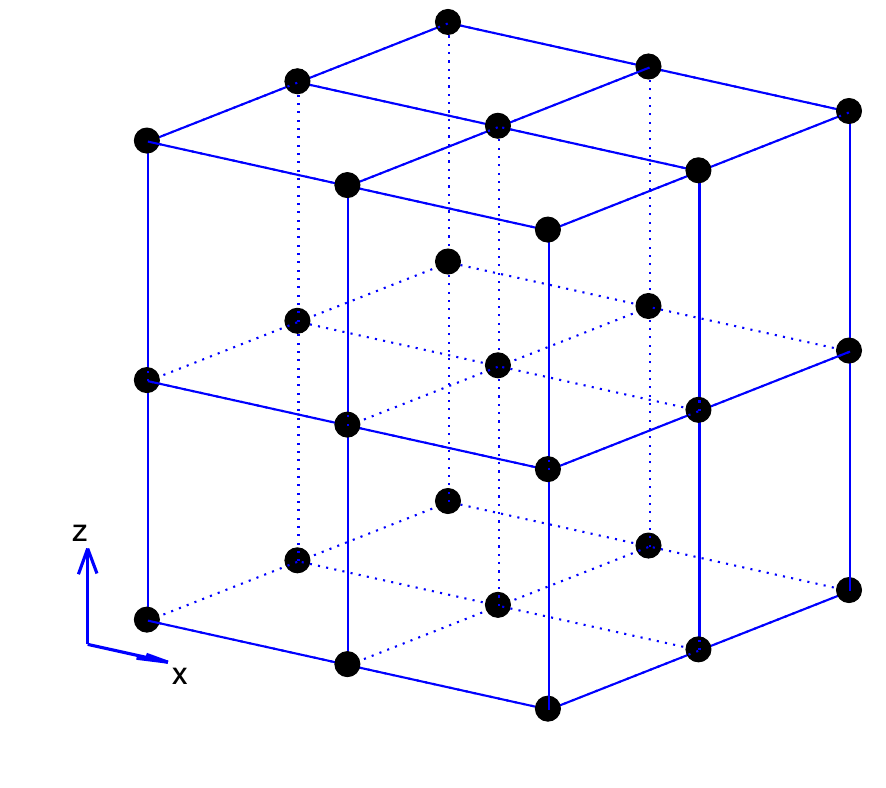}
}
\subfigure[$\mathcal{T}_h^0$]{
\includegraphics[width=0.25\textwidth]{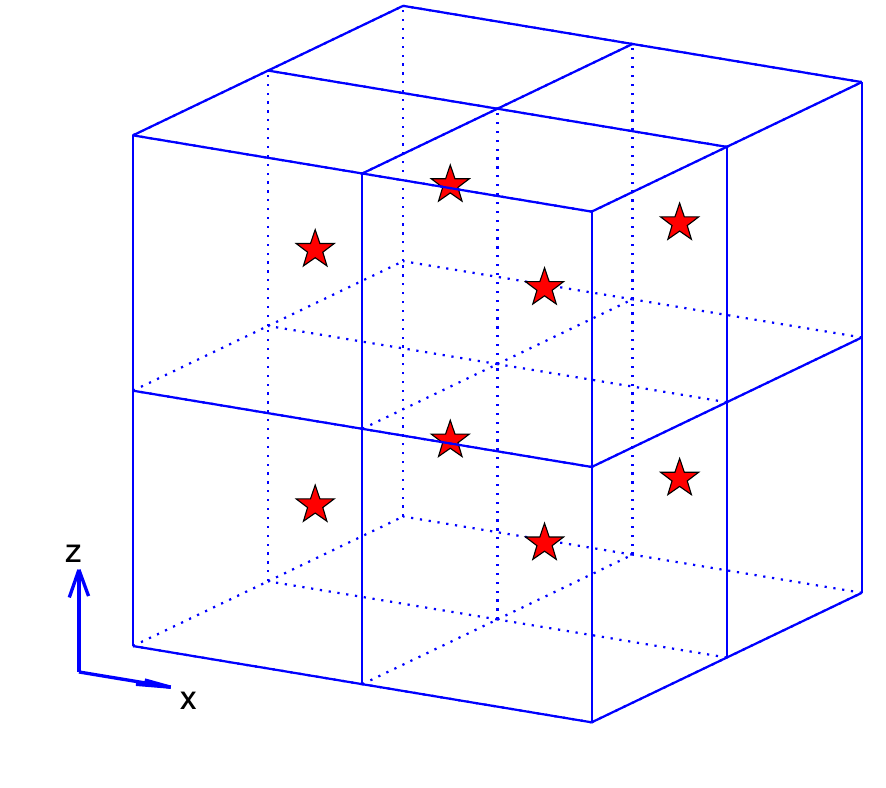}
}
\subfigure[$\mathcal{T}_h^1$]{
\includegraphics[width=0.25\textwidth]{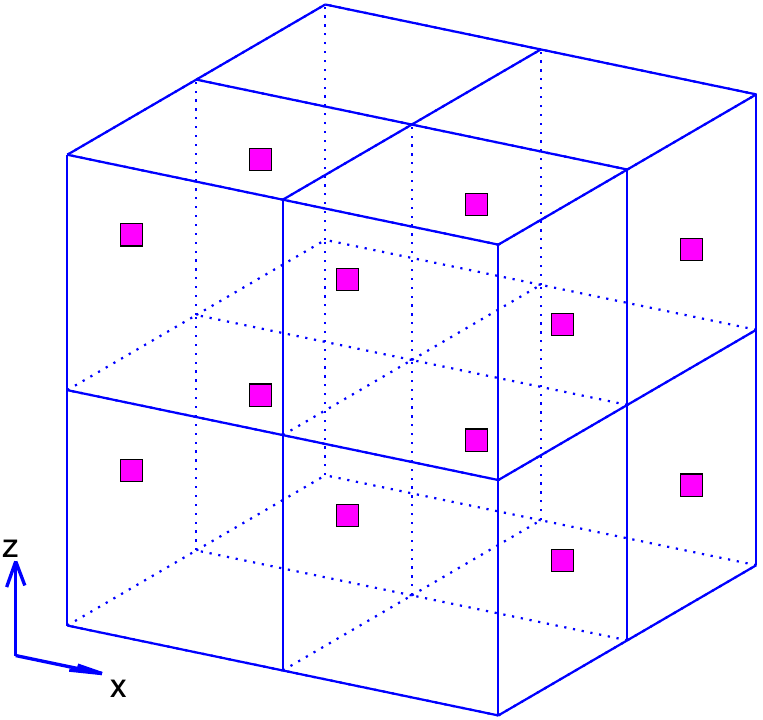}
}
\subfigure[$\mathcal{T}_h^2$]{
\includegraphics[width=0.25\textwidth]{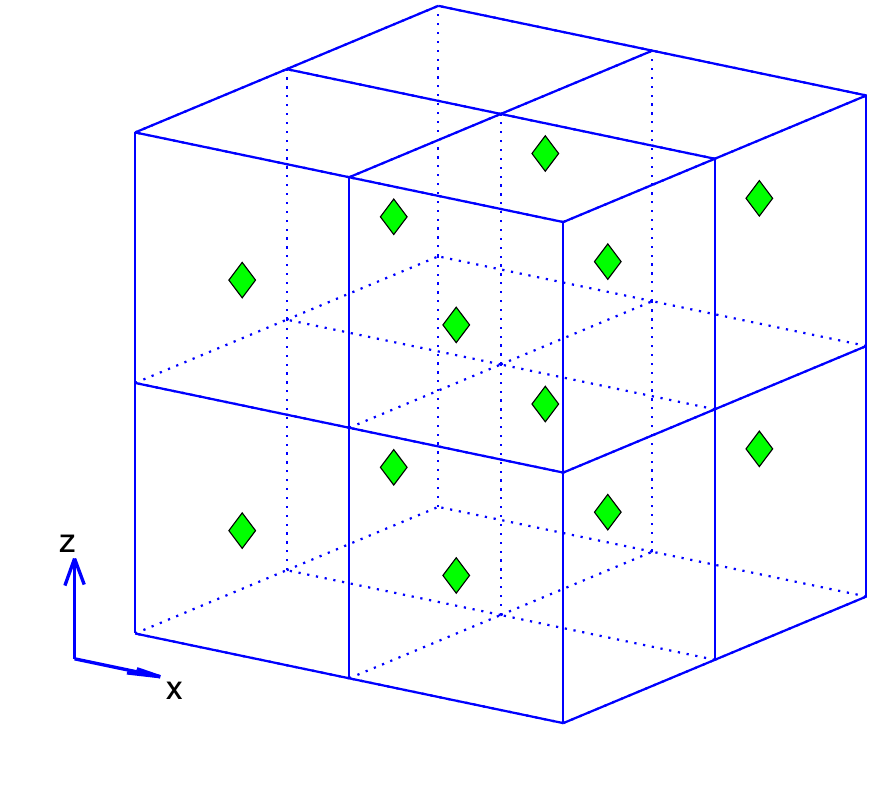}
}
\subfigure[$\mathcal{T}_h^3$]{
\includegraphics[width=0.25\textwidth]{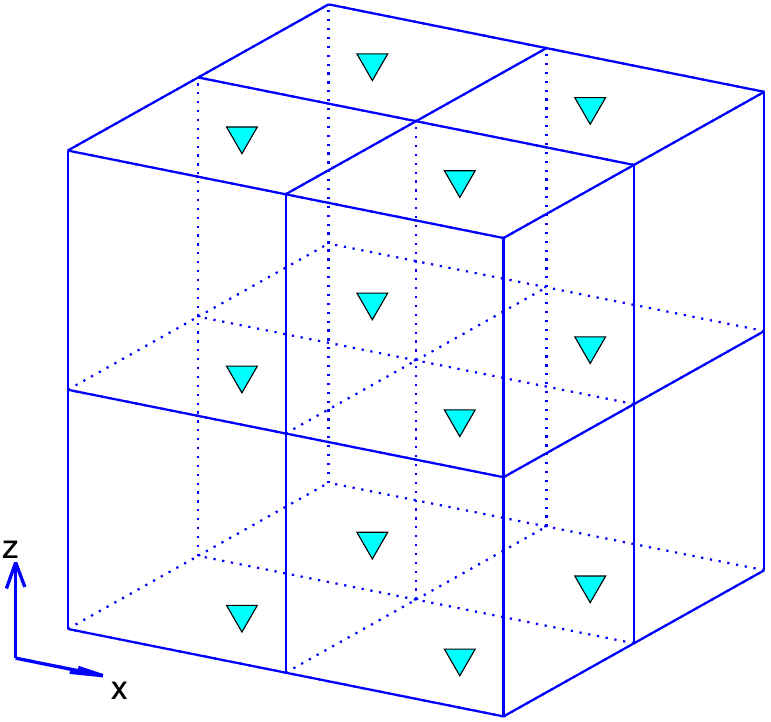}
}
\caption{Illustration of five different grid sets.}
\label{fig-mesh}
\end{figure}

Since there are no discontinuous coefficients in the interface problem \eqref{InterfaceUnfied}, numerical methods exist for solving it in the literature \cite{jones2021class, tan2011implementation, wang2021simple}. In this work,  a finite difference method on staggered grid with a delicate correction technique will be taken into consideration, which can be seen as a promotion of \cite{dong2022second} to 3D case. 
To simplify the presentation,  the subscript $`s$' in \eqref{InterfaceUnfied} will be omitted without confusion.
The computational domain $\mathcal{B}$ is taken as a unit cube, that is $\mathcal{B}\equiv [0,1]^3$. 
Suppose that the domain $\mathcal{B}$ is partitioned into $N\times N\times N$ uniform Cartesian grid with mesh parameter $h = x_{i+1} - x_i = y_{j+1}-y_j = z_{k+1}-z_k = 1/N $. For integers $i$, $j$, $k$, $0 \leq i \leq N, 0 \leq j \leq N, 0 \leq k \leq N$, define 
\begin{align*}
\begin{split}
x_{i-1/2} = \frac{x_i + x_{i-1}}{2}, \qquad y_{j-1/2}=\frac{y_j + y_{j-1}}{2}, \qquad z_{k-1/2}=\frac{z_k + z_{k-1}}{2}.
\end{split}
\end{align*}
Furthermore, introduce five different grid sets (see Fig. \ref{fig-mesh} for illustration): a vertex-centered grid set $\mathcal{T}_h$,  a cell-centered grid set $\mathcal{T}_h^0$, a YZ-plane-centered grid set $\mathcal{T}_h^1$,  a XZ-plane-centered grid set $\mathcal{T}_h^2$ and a XY-plane-centered grid set $\mathcal{T}_h^3$, given respectively by
\begin{align*}
	\begin{split}
\mathcal{T}_{h}&\equiv \big\{\,(x_i, y_j, z_k)\,|\,i=0, \cdots, N, \;j=0,  \cdots, N,\; k=0, \cdots, N \,\big\}, \\[4pt]
\mathcal{T}_{h}^0&\equiv \big\{\,(x_{i-\frac{1}{2}}, y_{j-\frac{1}{2}}, z_{k-\frac{1}{2}})\,|\,i=1,  \cdots, N, \;j=1,  \cdots, N,\;k=1,  \cdots, N\,\big\},\\[4pt]
\mathcal{T}_{h}^1&\equiv\big\{\,(x_i,y_ {j-\frac{1}{2}},z_{k-\frac{1}{2}})\,|\,i=0,  \cdots, N,\; j=0,  \cdots, N+1,\; k=0,  \cdots, N+1 \,\big\},\\[4pt]
\mathcal{T}_{h}^2&\equiv \big\{\,(x_{i-\frac{1}{2}},y_ j, z_{k-\frac{1}{2}})\,|\,i=0,  \cdots, N+1,\; j=0, \cdots, N,\; k=0,  \cdots, N+1\,\big\},\\[4pt]
\mathcal{T}_{h}^3&\equiv \big\{\,(x_{i-\frac{1}{2}},y_{j-\frac{1}{2}}, z_k)\,|\,i=0,  \cdots, N+1,\; j=0, \cdots, N+1,\; k=0,  \cdots, N\,\big\}.
	\end{split}
\end{align*}

A grid point $(x_{l},y_{m},z_{n})\in\mathcal{T}_h^{I},I=0,1,2,3,$ is called {\em irregular} point  if the corresponding finite difference stencils at this point go across the interface $\Gamma$, otherwise it is called {\em regular} point. Here $l,m$ and $n$ may take values $i, i-\frac{1}{2}, j,j-\frac{1}{2}$ and $k,k-\frac{1}{2}$.

For a function $v(x, y, z)$, set the differential operators
\begin{align*}
	\begin{split}
		\delta_{h,1}^{+}v_{l,m,n}=h^{-1}(v_{l+1,m,n}-v_{l,m,n}),\quad \quad \delta_{h,1}^{-}v_{l,m,n}=h^{-1}(v_{l,m,n}-v_{l-1,m,n}),\\[4pt]
		\delta_{h,2}^{+}v_{l,m,n}=h^{-1}(v_{l,m+1,n}-v_{l,m,n}),\quad \quad \delta_{h,2}^{-}v_{l,m,n}=h^{-1}(v_{l,m,n}-v_{l,m-1,n}),\\[4pt]
		\delta_{h,3}^{+}v_{l,m,n}=h^{-1}(v_{l,m,n+1}-v_{l,m,n}),\quad \quad \delta_{h,3}^{-}v_{l,m,n}=h^{-1}(v_{l,m,n}-v_{l,m,n-1}),
	\end{split}
\end{align*}
and the discrete Laplacian operator
\begin{align*}
	\begin{split}
		\Delta_h v_{l,m,n} = \delta_{h,1}^{+}\delta_{h,1}^{-}v_{l,m,n}+\delta_{h,2}^{+}\delta_{h,2}^{-}v_{l,m,n}+\delta_{h,3}^{+}\delta_{h,3}^{-}v_{l,m,n},
	\end{split}
\end{align*}
where $v_{l,m,n}$ represents $v(x_l, y_m, z_n)$.

Denote the exact solution of interface problem $\eqref{InterfaceUnfied}$ by
\begin{align*}
(\vect u,p) =& \Big(u^{(1)}(x_{i}, y_{j-\frac{1}{2}}, z_{k-\frac{1}{2}}), u^{(2)}(x_{i-\frac{1}{2}},y_{j}, z_{k-\frac{1}{2}}), u^{(3)}(x_{i-\frac{1}{2}}, y_{j-\frac{1}{2}}, z_{k}),\\
& p(x_{i-\frac{1}{2}}, y_{j-\frac{1}{2}}, z_{k-\frac{1}{2}})\Big),	
\end{align*}
and its finite difference approximation by 
$$(\vect u_h,p_h)=\big(u^{(1)}_{i,j-\frac{1}{2},k-\frac{1}{2}},
u^{(2)}_{i-\frac{1}{2}, j, k-\frac{1}{2}}, u^{(3)}_{i-\frac{1}{2}, j-\frac{1}{2},k}, p_{i-\frac{1}{2}, j-\frac{1}{2}, k-\frac{1}{2}}\big).$$ 
For simplicity of description, denote $I^{(1)}, \,I^{(2)},\, I^{(3)}, I^{(4)}$ as collections of indices for computational gird nodes, which are given by 
\begin{equation*}
\begin{split}
I^{(1)} &= \{(i,j,k)| \,  i=0,1,...,N -1, \; j= 1,2,...,N, \; k=1,2,...,N\},\\[2pt]
I^{(2)} &= \{(i,j,k)| \, i=1,2,...,N, \;  j= 0,1,...,N-1, \; k=1,2,...,N\},\\[2pt]
I^{(3)} &= \{(i,j,k)| \, i=1,2,...,N, \;  j= 1,2,...,N, \; k=0,1,...,N-1\},\\[2pt]
I^{(4)} &= \{(i,j,k)| \, i=1,2,...,N,\; j=1,2,...,N, \; k=1,2,...,N\}.
\end{split}
\end{equation*}
The discretization of the Stokes equations by second-order MAC scheme at a regular grid point reads as:
\begin{align*}
	\begin{split}
		-\mu\Delta_h u^{(1)}_{i,j-\frac{1}{2},k-\frac{1}{2}} + \delta^{+}_{h,1}\,p_{i-\frac{1}{2},j-\frac{1}{2},k-\frac{1}{2}}=\tilde{f}^{(1)}_{i,j-\frac{1}{2},k-\frac{1}{2}}, \quad\;\;(i,j,k)\in I^{(1)},\\[4pt]
		-\mu\Delta_h u^{(2)}_{i-\frac{1}{2},j,k-\frac{1}{2}} + \delta^{+}_{h,2}\,p_{i-\frac{1}{2},j-\frac{1}{2},k-\frac{1}{2}}=\tilde{f}^{(2)}_{i-\frac{1}{2},j,k-\frac{1}{2}},\quad\;\;(i,j,k)\in I^{(2)},\\[4pt]
		-\mu\Delta_h u^{(3)}_{i-\frac{1}{2},j-\frac{1}{2},k} + \delta^{+}_{h,3}\,p_{i-\frac{1}{2},j-\frac{1}{2},k-\frac{1}{2}}=\tilde{f}^{(3)}_{i-\frac{1}{2},j-\frac{1}{2},k},\quad\;\;(i,j,k)\in I^{(3)},\\[4pt]
		\delta^{-}_{h,1}u^{(1)}_{i,j-\frac{1}{2},k-\frac{1}{2}}+\delta^{-}_{h,2}u^{(2)}_{i-\frac{1}{2},j,k-\frac{1}{2}}+\delta^{-}_{h,3}u^{(3)}_{i-\frac{1}{2},j-\frac{1}{2},k}=0,\quad\;\;(i,j,k)\in I^{(4)}.
	\end{split}
\end{align*}
However, the above MAC scheme has large local truncation errors at an irregular point due to the existence of the interface $\Gamma$. To raise the global accuracy to second order,  appropriate correction terms  should be added to the right-hand side of the system, so that the coefficient matrix will not be changed. The modified  MAC scheme is in the form 
\begin{align}
\label{MAC}
	\begin{split}
		-\mu\Delta_h u^{(1)}_{i,j-\frac{1}{2},k-\frac{1}{2}} + \delta^{+}_{h,1}\,p_{i-\frac{1}{2},j-\frac{1}{2},k-\frac{1}{2}}=f^{(1)}_{i,j-\frac{1}{2},k-\frac{1}{2}}, \quad\;\;(i,j,k)\in I^{(1)},\\[4pt]
		-\mu\Delta_h u^{(2)}_{i-\frac{1}{2},j,k-\frac{1}{2}} + \delta^{+}_{h,2}\,p_{i-\frac{1}{2},j-\frac{1}{2},k-\frac{1}{2}}=f^{(2)}_{i-\frac{1}{2},j,k-\frac{1}{2}},\quad\;\;(i,j,k)\in I^{(2)},\\[4pt]
		-\mu\Delta_h u^{(3)}_{i-\frac{1}{2},j-\frac{1}{2},k} + \delta^{+}_{h,3}\,p_{i-\frac{1}{2},j-\frac{1}{2},k-\frac{1}{2}}=f^{(3)}_{i-\frac{1}{2},j-\frac{1}{2},k},\quad\;\;(i,j,k)\in I^{(3)},\\[4pt]
		\delta^{-}_{h,1}u^{(1)}_{i,j-\frac{1}{2},k-\frac{1}{2}}+\delta^{-}_{h,2}u^{(2)}_{i-\frac{1}{2},j,k-\frac{1}{2}}+\delta^{-}_{h,3}u^{(3)}_{i-\frac{1}{2},j-\frac{1}{2},k}=g_{i-\frac{1}{2},j-\frac{1}{2},k-\frac{1}{2}}\quad\;\;(i,j,k)\in I^{(4)},
	\end{split}
\end{align}
with 
\begin{align*}
	\begin{split}
		&f^{(1)}_{i,j-\frac{1}{2},k-\frac{1}{2}}=\tilde{f}^{(1)}_{i,j-\frac{1}{2},k-\frac{1}{2}}+C\{\mu\Delta u^{(1)}\}_{i,j-\frac{1}{2},k-\frac{1}{2}}+C\{p_x\}_{i,j-\frac{1}{2},k-\frac{1}{2}},\\[4pt]
		&f^{(2)}_{i-\frac{1}{2},j,k-\frac{1}{2}}=\tilde{f}^{(2)}_{i-\frac{1}{2},j,k-\frac{1}{2}}+C\{\mu\Delta u^{(2)}\}_{i-\frac{1}{2},j,k-\frac{1}{2}}+C\{p_y\}_{i-\frac{1}{2},j,k-\frac{1}{2}},\\[4pt]
		&f^{(3)}_{i-\frac{1}{2},j-\frac{1}{2},k}=\tilde{f}^{(3)}_{i-\frac{1}{2},j-\frac{1}{2},k}+C\{\mu\Delta u^{(3)}\}_{i-\frac{1}{2},j-\frac{1}{2},k}+C\{p_z\}_{i-\frac{1}{2},j-\frac{1}{2},k},\\[4pt]
		&g_{i-\frac{1}{2},j-\frac{1}{2},k-\frac{1}{2}}=C\{u^{(1)}_x\}_{i-\frac{1}{2},j-\frac{1}{2},k-\frac{1}{2}}+C\{u^{(2)}_y\}_{i-\frac{1}{2},j-\frac{1}{2},k-\frac{1}{2}}+C\{u^{(3)}_z\}_{i-\frac{1}{2},j-\frac{1}{2},k-\frac{1}{2}}.
	\end{split}
\end{align*}
It is noted that the correction terms are the sum of a few leading order terms of Taylor expansions and are non-zero only at irregular points. And they will improve the local truncation errors near the interface to at least first-order accuracy. As one will see in the coming section, these correction terms can be computed in terms of the jumps of the solution and their partial derivatives. As a matter of fact, all jump conditions are also computable by taking the derivatives of the interface conditions together with the incompressible condition. 

Let $\pmb\Delta_h={\rm diag}(\Delta_h, \Delta_h, \Delta_h)$ denote the standard central difference operator, and $G_h^{\rm MAC}=(\delta^{+}_{h,1},\delta^{+}_{h,2},\delta^{+}_{h,3})^T,$ $ D_h^{\rm MAC}=(\delta^{-}_{h,1},\delta^{-}_{h,2},\delta^{-}_{h,3})$ be the MAC gradient or divergence operator, then the MAC scheme can be rewritten as a linear system in the form of 
\begin{equation}
\label{linearS}
\begin{pmatrix}
-\mu\pmb \Delta_h  & G_h^{\rm MAC}  \\[4pt]
D_h^{\rm MAC} &  \vect 0\\
\end{pmatrix}
\begin{pmatrix}
\vect u_h\\[4pt]
p_h\\
\end{pmatrix}=
\begin{pmatrix}
\vect f\;\\[4pt]
g\\
\end{pmatrix}.
\end{equation}
It is known that the resulting system is a saddle point problem, which must be solved using an iterative method, typically a Krylov subspace method such as GMRES \cite{saad1986gmres}.   A wide variety of preconditioners have been proposed for such systems, mainly including domain decomposition methods \cite{lin2006performance,cyr2012stabilization},  block preconditioners \cite{elman2006block,griffith2009accurate,bootland2019preconditioners} and  multigrid methods \cite{vanka1986block,elman1996multigrid,wang2013multigrid}.  Here, a conjugate gradient (CG) method together with an FFT-based Poisson solver introduced in our previous work \cite{dong2022second} has been extended to the present three-dimensional case, which can be included as
\begin{compactitem}
\item[i)]Since the trouble is the uniqueness of the pressure $p_h$, an auxiliary variable $\lambda_h$ and a parameter $\alpha$ satisfying the condition that $\lambda_h$ equals the average of the pressure variable over the domain are introduced to ensure the solvability of the system. Then the linear system reads as
\begin{equation}
\label{linearS1}
\begin{pmatrix}
-\mu\pmb \Delta_h  & G_h^{\rm MAC} & 0 \\[4pt]
D_h^{\rm MAC} &  0 & -\gamma \\[4pt]
0 & -\gamma^T & \alpha
\end{pmatrix}
\begin{pmatrix}
\vect u_h\\[4pt]
p_h\\[4pt]
\lambda_h
\end{pmatrix}=
\begin{pmatrix}
\vect f\\[4pt]
g\\[4pt]
0
\end{pmatrix}.
\end{equation}
\item[ii)] Substituting $\vect u_h = (-\Delta_h^{-1}\vect f+\Delta_h^{-1}G_h^{\rm MAC}p_h)/\mu$ into \eqref{linearS1} yields the system
\begin{equation*}
(\frac{1}{\mu}D_h^{\rm MAC}\pmb \Delta_h^{-1}G_h^{\rm MAC}+\dfrac{1}{\alpha}\gamma\gamma^T)p_h=\frac{1}{\mu}D_h^{\rm MAC}\Delta_h^{-1}\vect f+g,
\end{equation*}
which can be solved efficiently by the CG method. Besides, the evaluation of $\Delta_h^{-1}$ can be transformed into solving a Poisson equation. Hence an FFT-based Poisson fast solver can be used in each iteration. 
\item[iii)] Furthermore, using $p_h$ computed from the above system, the velocity $\vect u_h$ can be solved from
 \begin{equation*}
 -\mu\vect \Delta_h \vect u_h=\vect f-G_h^{\rm MAC}p_h,
 \end{equation*}
with the FFT-based Poisson solvers again.
\end{compactitem}

It is remarked that no preconditioner for the iterative method is involved in the present method. Efficient algorithms incorporated with preconditioner for the time-dependent problems are encouraging and will be reported in future work.

\subsection{Correction of the MAC Scheme}

Once again, the discontinuity of the velocity and its stress tensor across the interface $\Gamma$ leads to the fact that the local truncation errors of the finite difference MAC scheme \eqref{MAC} at irregular grid nodes are too large, and the solution to the discrete Stokes interface problem will be inaccurate without any modification. To obtain the desired second-order accuracy, some correction terms have been added to the right hand of the discrete system \eqref{MAC}.  This subsection presents the detailed derivation of the correction terms. Here, the correction technique is similar to that presented in \cite{ying2013kernel} except that the used grid is staggered grid and the pressure term $p$ is also needed to be modified. For simplicity, only the right side of the x-axis is illustrated, and the correction terms at the opposite sides can be obtained by symmetry. Moreover, the correction terms at the y- or z-direction can be obtained using the same method.  Assuming that there is an irregular point on the gird cell in Fig.\ref{fig-irregular} (a), the correction terms are evaluated as follows:
\begin{figure}[!ht]
\centering
\subfigure[Reference Cell]{
\includegraphics[width=0.25\textwidth]{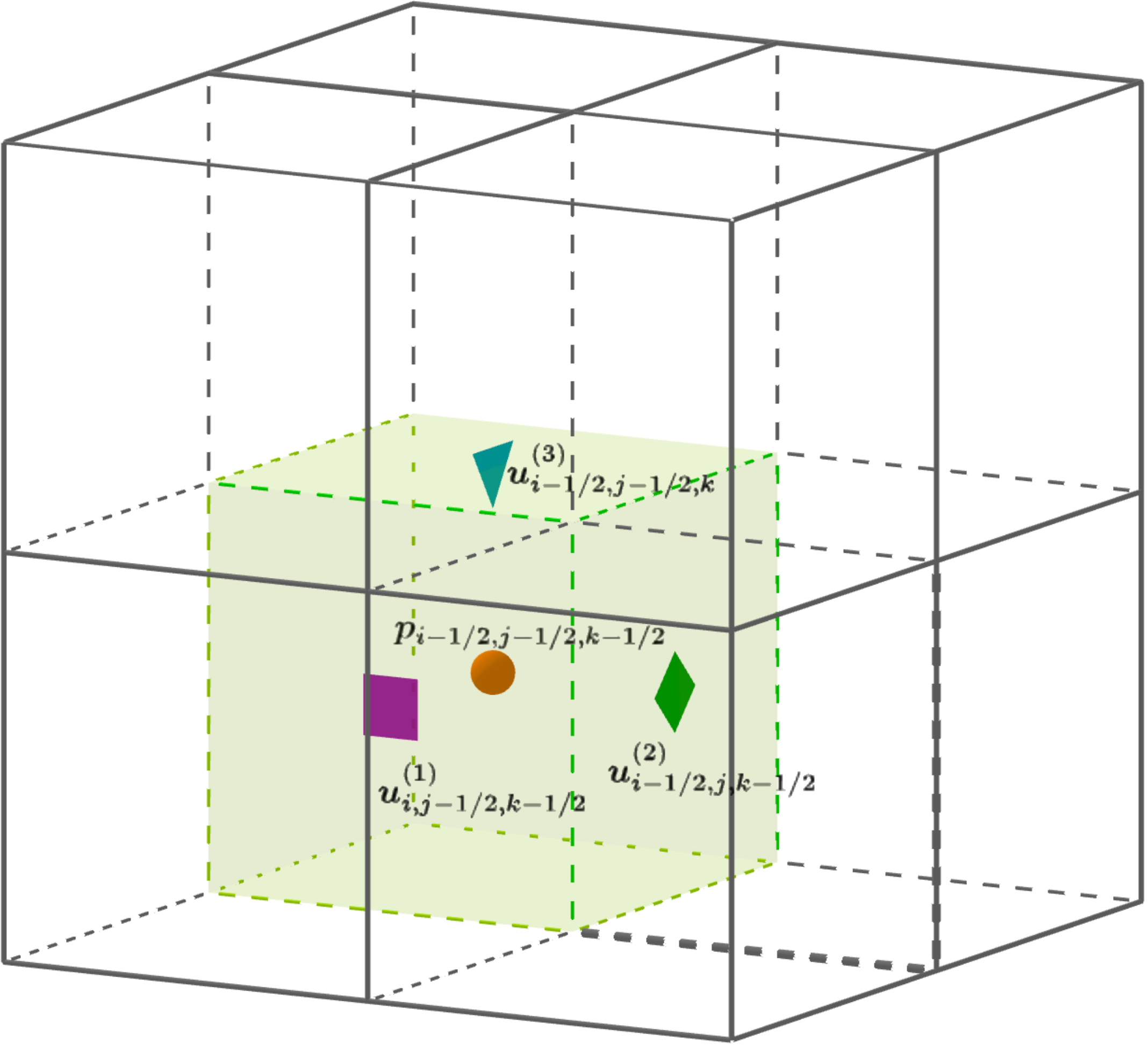}
}
\subfigure[$(x_i,y_{j-\frac{1}{2}},z_{k-\frac{1}{2}})$]{
\includegraphics[width=0.25\textwidth]{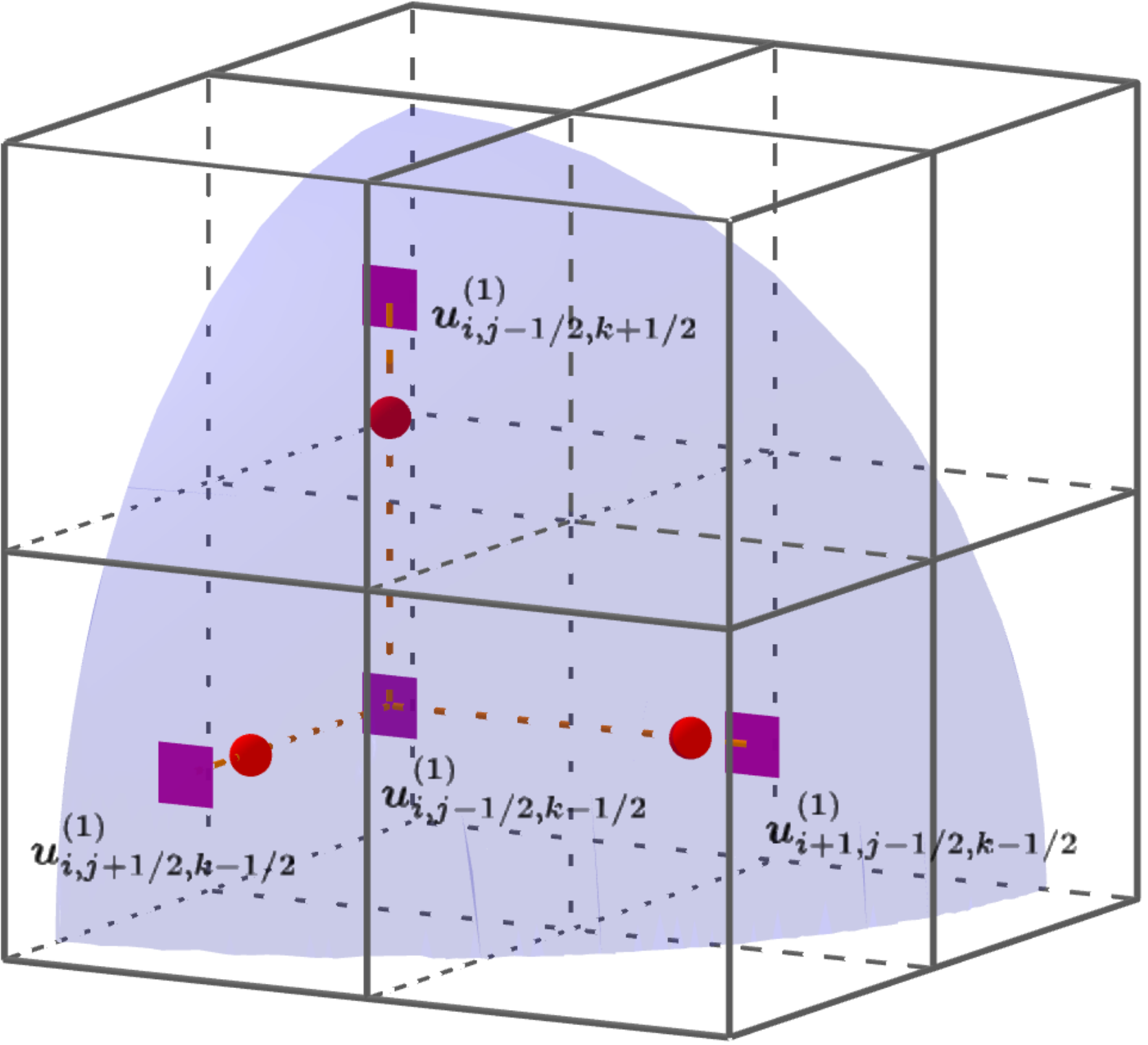}
}
\subfigure[$(x_{i-\frac{1}{2}},y_j,z_{k-\frac{1}{2}})$]{
\includegraphics[width=0.25\textwidth]{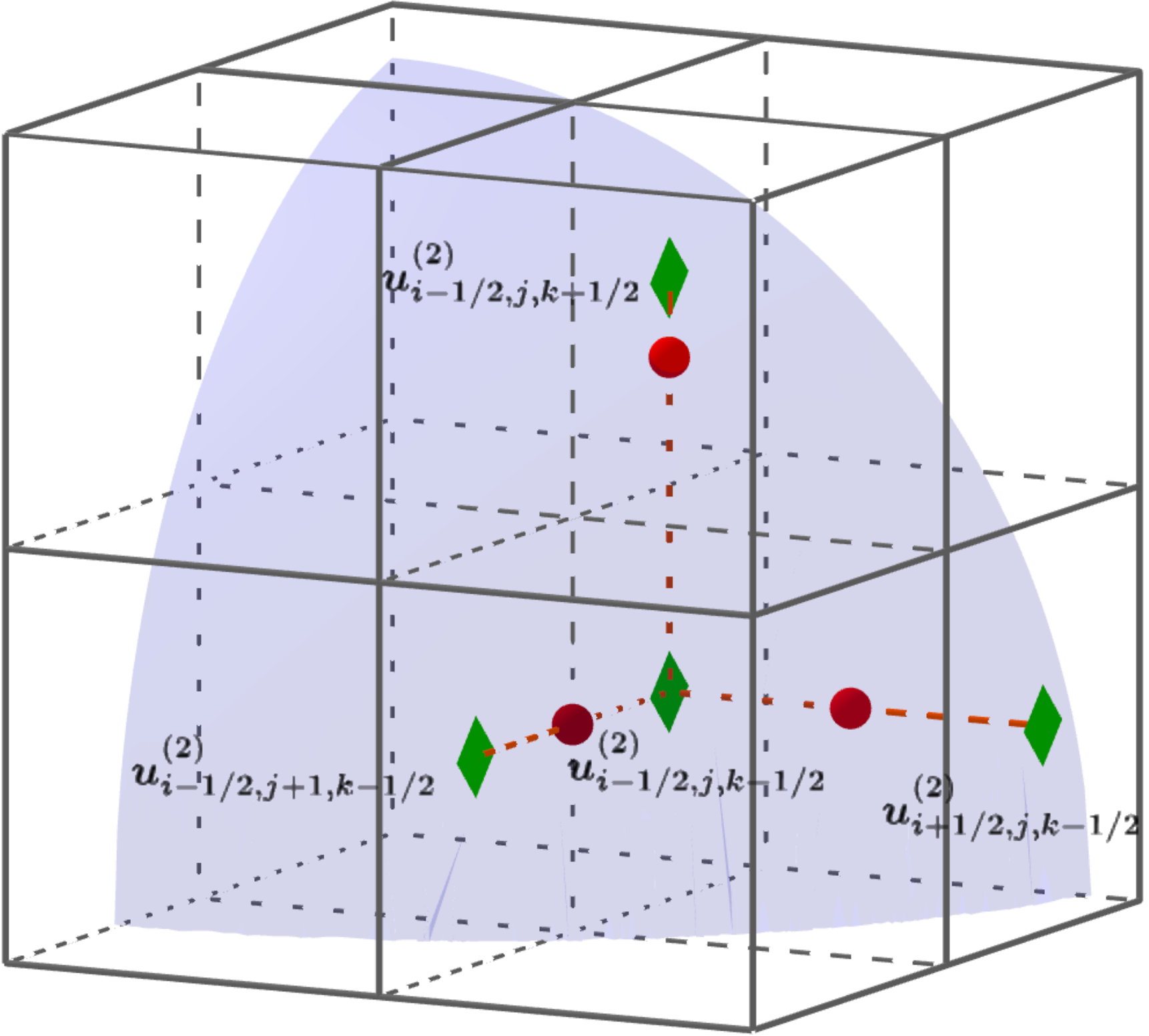}
}
\subfigure[$(x_{i-\frac{1}{2}},y_{j-\frac{1}{2}},z_k)$]{
\includegraphics[width=0.25\textwidth]{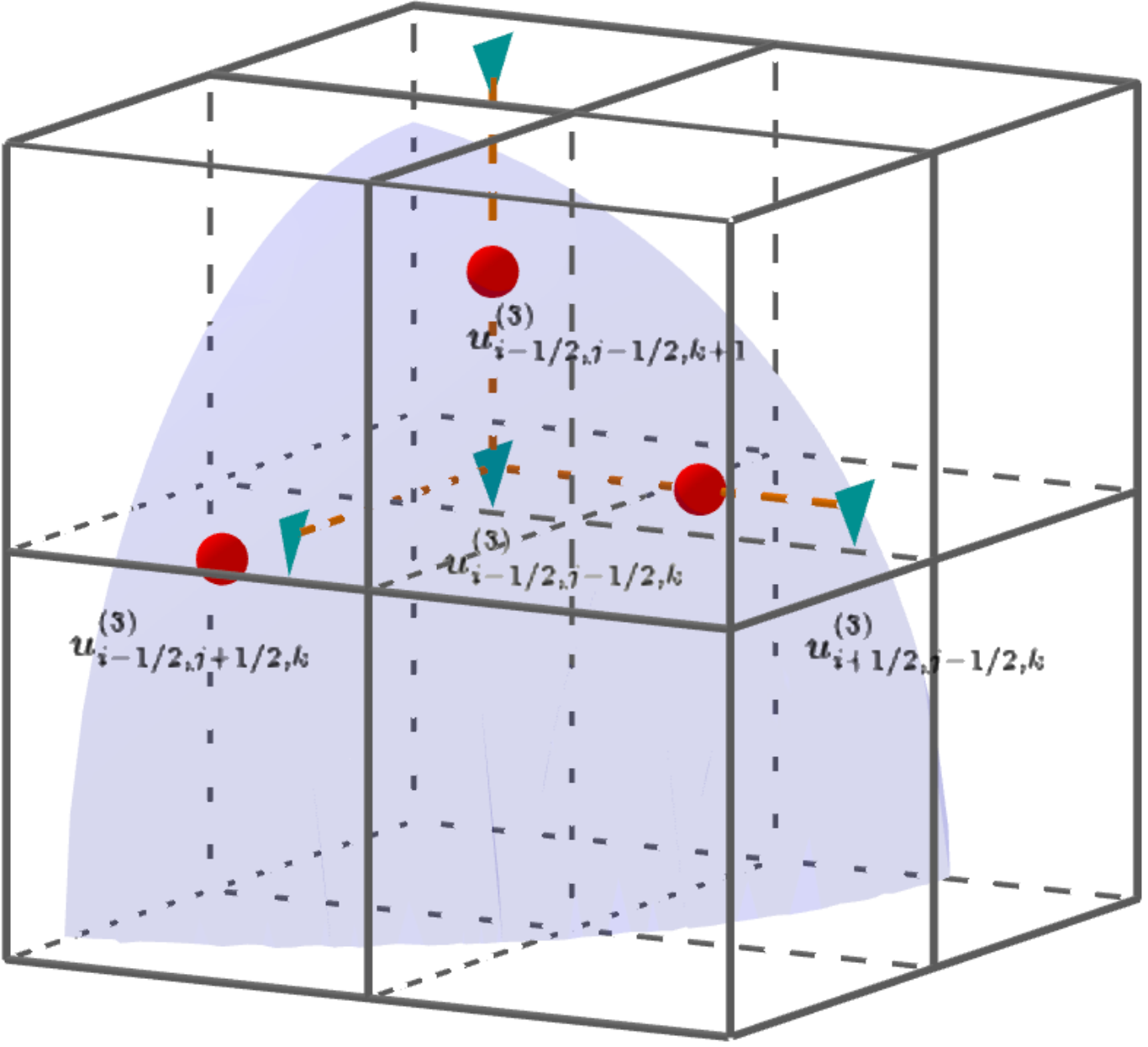}
}
\subfigure[$(x_{i-\frac{1}{2}},y_{j-\frac{1}{2}},z_{k-\frac{1}{2}})$]{
\includegraphics[width=0.25\textwidth]{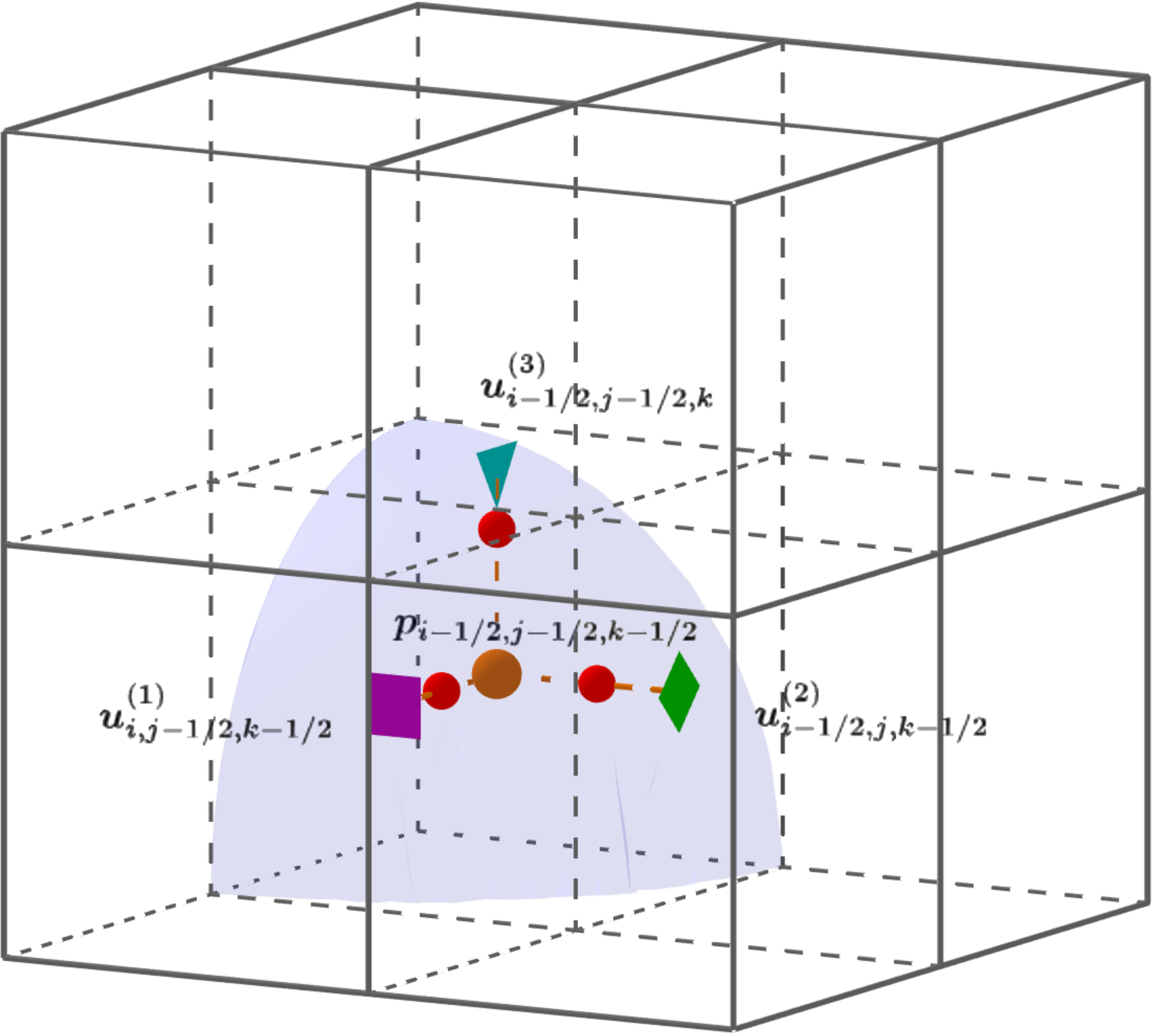}
}
\caption{Illustration of Irregular points.}
\label{fig-irregular}
\end{figure}


\noindent 1. {\em Assume that $(x_i,y_{j-\frac{1}{2}},z_{k-\frac{1}{2}})$ is an irregular point, see Fig.$\ref{fig-irregular} (b)$ for illustration.} In the $x$-direction, denote by $\xi_{u^{(1)}}=x_{i+1}-r$. If $\Gamma$ intersects the grid line segment consisting of point $(x_i, y_{j-\frac{1}{2}}, z_{k-\frac{1}{2}})$ and $(x_{i+1}, y_{j-\frac{1}{2}}, z_{k-\frac{1}{2}})$ at $(r, y_{j-\frac{1}{2}}, z_{k-\frac{1}{2}})$ with $x_{i} < r < x_{i+1}$, the correction term can be computed from Taylor expansion  as
\begin{equation*}
C\{\mu\Delta u^{(1)}\}_{i,j-\frac{1}{2},k-\frac{1}{2}} = 
\dfrac{\mu}{h^2}\Big( [\![ u^{(1)} ]\!]  
+ \xi_{u^{(1)}} [\![ u^{(1)}_x ]\!] 
+\dfrac{1}{2}\xi_{u^{(1)}}^2 [\![ u^{(1)}_{xx} ]\!]   \Big), 
\ \  \hbox{if}\; (x_i,y_{j-\frac{1}{2}},z_{k-\frac{1}{2}})\in \Omega.
\end{equation*}

\noindent2. {\em Assume that $(x_{i-\frac{1}{2}},y_j,z_{k-\frac{1}{2}})$ is an irregular point, see Fig.$\ref{fig-irregular}(c)$ for illustration.} In the $x$-direction, denote by $\xi_{u^{(2)}}=x_{i+\frac{1}{2}}-r$. If $\Gamma$ intersects the grid line segment consisting of point $(x_{i-\frac{1}{2}},y_j,z_{k-\frac{1}{2}})$ and $(x_{i+\frac{1}{2}},y_j,z_{k-\frac{1}{2}})$ at $(r,y_j,z_{k-\frac{1}{2}})$ with $x_{i-\frac{1}{2}}<r< x_{i+\frac{1}{2}}$, the correction term can  be  
computed   from Taylor expansion  as
\begin{equation*}
C\{\mu\Delta u^{(2)}\}_{i-\frac{1}{2},j,k-\frac{1}{2}} = 
\dfrac{\mu}{h^2}\Big( [\![ u^{(2)} ]\!]  
+ \xi_{u^{(2)}} [\![ u^{(2)}_x ]\!] 
+\dfrac{1}{2}\xi_{u^{(2)}}^2 [\![ u^{(2)}_{xx} ]\!]   \Big), 
\ \ \hbox{if}\; (x_{i-\frac{1}{2}},y_j,z_{k-\frac{1}{2}})\in \Omega.
\end{equation*}

\noindent 3. {\em Assume that $(x_{i-\frac{1}{2}},y_{j-\frac{1}{2}},z_k)$ is an irregular point, see Fig. $\ref{fig-irregular} (d)$ for illustration.} In the $x$-direction, denote by $\xi_{u^{(3)}}=x_{i+\frac{1}{2}}-r$. $\Gamma$ intersects the grid line segment consisting of point $(x_{i-\frac{1}{2}},y_{j-\frac{1}{2}},z_k)$ and $(x_{i+\frac{1}{2}},y_{j-\frac{1}{2}},z_k)$ at $(r,y_{j-\frac{1}{2}},z_k)$ with $x_{i-\frac{1}{2}} < r < x_{i+\frac{1}{2}}$, the correction term can be computed from Taylor expansion as
\begin{equation*}
C\{\mu\Delta u^{(3)}\}_{i-\frac{1}{2},j-\frac{1}{2},k} = 
\dfrac{\mu}{h^2}\Big( [\![ u^{(3)} ]\!]  
+ \xi_{u^{(3)}} [\![ u^{(3)}_x ]\!] 
+\dfrac{1}{2}\xi_{u^{(3)}}^2 [\![ u^{(3)}_{xx} ]\!]   \Big), 
\ \  \hbox{if}\; (x_{i-\frac{1}{2}},y_{j-\frac{1}{2}},z_k)\in \Omega.
\end{equation*}

\noindent 4. {\em Assume that $(x_{i-\frac{1}{2}},y_{j-\frac{1}{2}},z_{k-\frac{1}{2}})$ is an irregular point, see Fig.$\ref{fig-irregular} (e)$  for illustration.}
 In the $x$-direction, denote by $\xi_{u^{(1)}}=x_i-r$ and $\xi_p = x_{i-\frac{1}{2}} - r$. If $\Gamma$ intersects the grid line segment consisting of point $(x_{i-\frac{1}{2}},y_{j-\frac{1}{2}},z_{k-\frac{1}{2}})$ and $(x_i,y_{j-\frac{1}{2}},z_{k-\frac{1}{2}})$ at $(r,y_{j-\frac{1}{2}},z_{k-\frac{1}{2}})$ with $x_{i-\frac{1}{2}} < r < x_i$, the correction term can be computed   from Taylor expansion as
\begin{align*}
C\{u^{(1)}_x\}_{i-\frac{1}{2},j-\frac{1}{2},k-\frac{1}{2}} = 
\dfrac{1}{h}\Big( [\![ u^{(1)} ]\!]  
+ \xi_{u^{(1)}} [\![ u^{(1)}_x ]\!] 
&+\dfrac{1}{2}\xi_{u^{(1)}}^2 [\![ u^{(1)}_{xx} ]\!]   \Big), \\
&\hbox{if}\; (x_{i-\frac{1}{2}},y_{j-\frac{1}{2}},z_{k-\frac{1}{2}})\in \Omega,
\end{align*}
{\em and} 
\begin{equation*}
C\{p_x\}_{i,j-\frac{1}{2},k-\frac{1}{2}} = 
\dfrac{1}{h}\Big( [\![ p ]\!]  
+ \xi_{p} [\![ p_x ]\!]  \Big), 
\quad \hbox{if}\; (x_i,y_{j-\frac{1}{2}},z_{k-\frac{1}{2}})\in \Omega.
\end{equation*}

Once the jumps across the interface are computed, the correction terms can be derived explicitly, so that they can be added to the right hand of the discrete linear system \eqref{MAC} at the irregular grid nodes. Thus, the coefficient matrix of the modified system is the same as the standard Stokes problem without an interface and the existing fast solver is still applicable. Besides, the derivation of the correction terms indicates that they improved the local truncation errors near the interface to at least first order accuracy, which is sufficient
 to recover the formal second-order accuracy of the underlying numerical scheme. Numerical results in section \ref{sec;numerical} illustrate this fact, and one can refer \cite{dong2022second} for the detailed analysis.
 
In addition, the jump conditions involved in the correction terms can be calculated from the original equations \eqref{IUS-1}-\eqref{IUS-4}, and the details will not be presented here because of the size limitation.

\subsection{Interpolation for Integrals on the Interface}
\begin{figure}[!ht]
  \centering
    \includegraphics[width=0.4\textwidth]{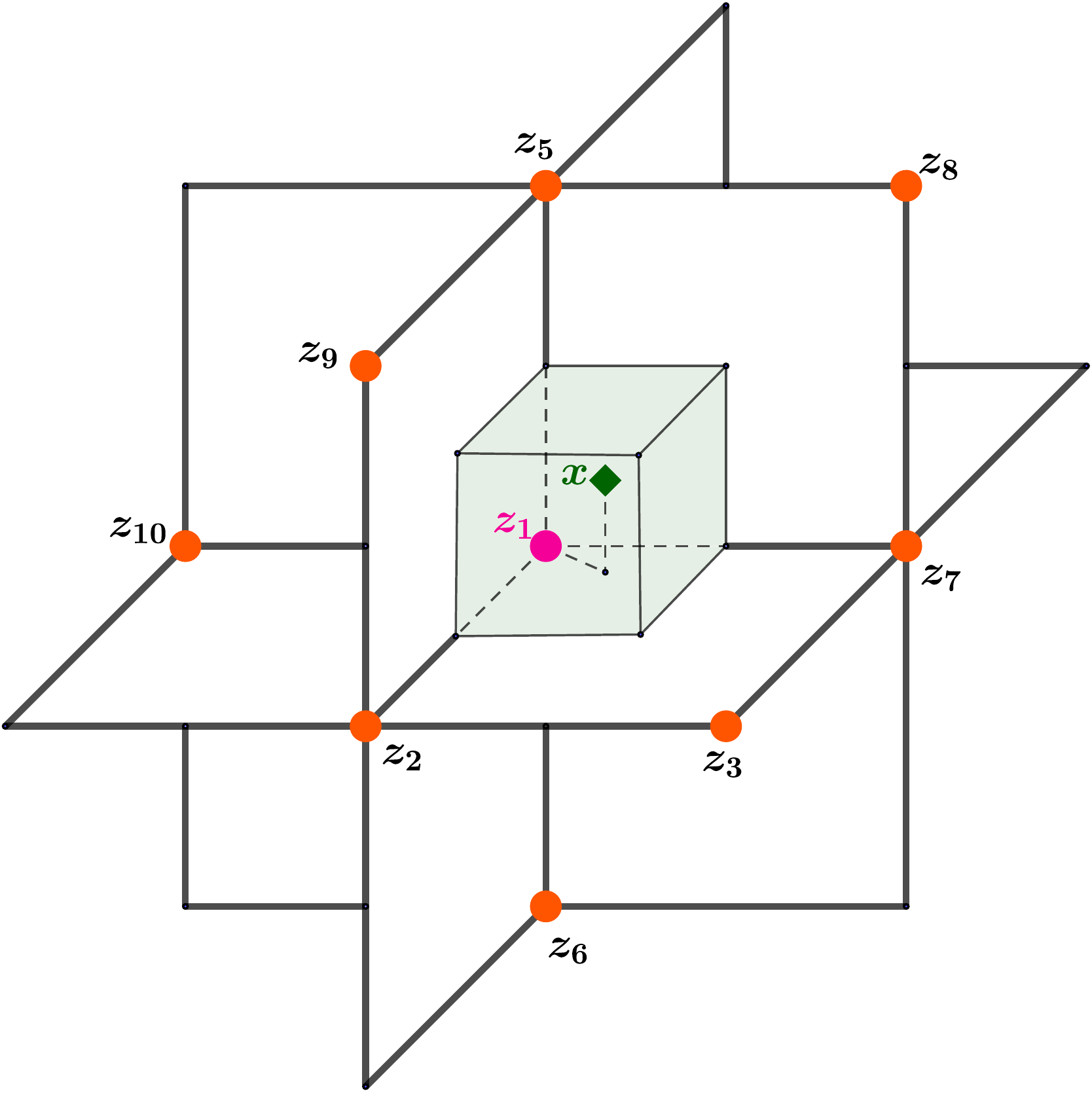}
\caption{
Point $\vect x$ on the interface located at one quadrant distribution and the corresponding ten interpolation grid nodes $\vect z_i$ (nine of them are marked with circle and the other one ${\vect z}_4$ is hidden behind) for computing the limits values of an approximation solution and its traction.}
\label{fig-interpolation}
\end{figure}

It is seen that the approximation solution to the interface problem \eqref{InterfaceUnfied} is calculated on a staggered grid node, while the approximations of the corresponding boundary or volume integral needed in \eqref{BIE-S} should be evaluated at discretization points of the interface $\Gamma$. To this end,  with $(\vect u_h, p_h)$  computed by the MAC scheme \eqref{MAC},  a quadratic polynomial  interpolation should be designed to extract limit values of $\vect u_h$ and its flux $\pmb\sigma(\vect u,p)$ at any given discretization points $\vect x$ on the interface. 
For a control point $\vect x $ on the interface, ten closest grid nodes $\vect z_i (i = 0,1,2, ... ,9)$ are chosen to construct the interpolation stencil (see Fig. \ref{fig-interpolation} for illustration). 

For each interpolation point $\vect z_i \in \Omega (i = 1,2, ... ,10)$, Taylor expansion around the point $\vect x\in\Gamma$ is given by
\begin{align}
\label{Interpolation+}
\begin{split}
	\vect u_h(\vect z_i) =& \vect u_h^+(\vect x)+  \xi_i\vect u_{h,x}^+(\vect x)+ \eta_i\vect u_{h,y}^+(\vect x)+  \gamma_i\vect u_{h,z}^+(\vect x) \\[4pt]
	&+ \frac{1}{2}\xi_i^2\vect u_{h,xx}^+(\vect x)
	+ \frac{1}{2}\eta_i^2\vect u_{h,yy}^+(\vect x)
	+ \frac{1}{2}\gamma_i^2\vect u_{h,zz}^+(\vect x)\\[4pt]
	&+ \eta_i\gamma_i \vect u_{h,yz}^+(\vect x)
	+ \xi_i\gamma_i \vect u_{h,xz}^+(\vect x)
	+ \xi_i\eta_i \vect u_{h,xy}^+(\vect x)+O(|\vect z_i - \vect x|^3),
\end{split}
\end{align}
where $\vect z_i - \vect x \equiv (\xi_i, \eta_i, \gamma_i)^T$, $\vect u_{h,x}$, $\vect u_{h,y}$, $\vect u_{h,z}$ are respectively used to denote $\partial \vect u_h/\partial x$, $\partial \vect u_h/\partial y$ and $\partial \vect u_h/\partial z$ and the second order derivatives are defined in a similar way.  It is noted that $\vect z_i$ belongs to different grid sets $\mathcal{T}_h^{j} (j=1,2,3)$ for different components $u^{(j)}$.
While for each interpolation point $\vect z_i \in \Omega^c$, Taylor expansion around $\vect x$ is given by 
\begin{align}
\label{Interpolation-}
\begin{split}
	\vect u_h(\vect z_i) =& \vect u_h^-(\vect x)+  \xi_i\vect u_{h,x}^-(\vect x)+ \eta_i\vect u_{h,y}^-(\vect x)+  \gamma_i\vect u_{h,z}^-(\vect x) \\[4pt]
	&+ \frac{1}{2}\xi_i^2\vect u_{h,xx}^-(\vect x)
	+ \frac{1}{2}\eta_i^2\vect u_{h,yy}^-(\vect x)
	+ \frac{1}{2}\gamma_i^2\vect u_{h,zz}^-(\vect x)\\[4pt]
	&+ \eta_i\gamma_i \vect u_{h,yz}^-(\vect x)
	+ \xi_i\gamma_i \vect u_{h,xz}^-(\vect x)
	+ \xi_i\eta_i \vect u_{h,xy}^-(\vect x)+O(|\vect z_i - \vect x|^3),
\end{split}
\end{align}
For conciseness, denote the approximate value $\vect u_{h}$ and its derivatives respectively by
\begin{align*}
	\begin{split}
		&\vect v^{\pm}=\vect u^{\pm}_{h}(\vect x), \ \ \ \qquad
		\vect v_x^{\pm}=\vect u^{\pm}_{h,x}(\vect x), \ \ \qquad
		\vect v_y^{\pm}=\vect u^{\pm}_{h,y}(\vect x), \\[4pt]
		&\vect v_{xx}^{\pm}=\vect u^{\pm}_{h,xx}(\vect x),  \qquad
		\vect v_{yy}^{\pm}=\vect u^{\pm}_{h,yy}(\vect x), \qquad
		\vect v_{zz}^{\pm}=\vect u^{\pm}_{h,zz}(\vect x), \\[4pt]
		&\vect v_{yz}^{\pm}=\vect u^{\pm}_{h,yz}(\vect x),  \qquad
		\vect v_{xz}^{\pm}=\vect u^{\pm}_{h,xz}(\vect x), \qquad
		\vect v_{xy}^{\pm}=\vect u^{\pm}_{h,xy}(\vect x).
	\end{split}
\end{align*}

Due to the discontinuity of $\vect v(\vect x)$ and its traction across the interface $\Gamma$, the following jump function should be introduced 
\begin{align*}
\begin{split}
	\vect J_i =& [\![ \vect v]\!] + \xi_i[\![ \vect v_x]\!] + \eta_i [\![ \vect v_y]\!] + \gamma_i[\![ \vect v_z]\!] 
	+\frac{1}{2}\xi_i^2[\![ \vect v_{xx}]\!] 
	+\frac{1}{2}\eta_i^2[\![ \vect v_{yy}]\!] 
	+\frac{1}{2}\gamma_i^2[\![ \vect v_{zz}]\!] \\[4pt]
	&+\eta_i\gamma_i[\![ \vect v_{yz}]\!] 
	+\xi_i\gamma_i[\![ \vect v_{xz}]\!] 
	+\xi_i\eta_i[\![ \vect v_{xy}]\!].
\end{split}
\end{align*}
Thus, Taylor expansion (\ref{Interpolation+}) and (\ref{Interpolation-}) can be rewritten as 
\begin{align*}
\begin{split}
	\;\qquad\vect v_i =& \vect v^{+} +  \xi_i\vect v_{x}^{+}+ \eta_i\vect v_{y}^{+}+  \gamma_i\vect v_{z}^{+} 
	+ \frac{1}{2}\xi_i^2\vect v_{xx}^{+}
	+ \frac{1}{2}\eta_i^2\vect v_{yy}^{+}
	+ \frac{1}{2}\gamma_i^2\vect v_{zz}^{+}\\[4pt]
    &+ \eta_i\gamma_i \vect v_{yz}^{+}
	+ \xi_i\gamma_i \vect v_{xz}^{+}
	+ \xi_i\eta_i \vect v_{xy}^{+},\quad {\rm if} \ \vect z_i \in \Omega,
\end{split}
\end{align*}
and
\begin{align*}
\begin{split}
	\vect v_i + \vect J_i =& \vect v^{+} +  \xi_i\vect v_{x}^{+}+ \eta_i\vect v_{y}^{+}+  \gamma_i\vect v_{z}^{+} 
	+ \frac{1}{2}\xi_i^2\vect v_{xx}^{+}
	+ \frac{1}{2}\eta_i^2\vect v_{yy}^{+}
	+ \frac{1}{2}\gamma_i^2\vect v_{zz}^{+}\\[4pt]
   &+ \eta_i\gamma_i \vect v_{yz}^{+}
	+ \xi_i\gamma_i \vect v_{xz}^{+}
	+ \xi_i\eta_i \vect v_{xy}^{+},\quad {\rm if} \ \vect z_i \in \Omega^c.
\end{split}
\end{align*}
Here, the third-order term $O(|\vect z_i - \vect x|^3)$ is omitted.

Solving the linear system above, one could obtain the limit value of $\vect u_h$ on the interface $\Gamma$. Since the stress tensor $\pmb \sigma(\vect u_h, p_h)$ does not involve the derivatives of $p_h$, linear interpolation is enough in the evaluation of boundary value $p_h$, which can be done similarly and will be omitted here.

\section{Numerical Results}
\label{sec;numerical}
Numerical tests are designed in this section to investigate the accuracy, efficiency and robustness of the proposed method for solving 3D Stokes and Navier boundary value problems. To this end, the normalized maximum norms and $\ell^2$-norms are defined as follows 
$$
\|e_{p}\|_{\infty}=\frac{\|p-p_h\|_{\infty}}{\|p\|_{\infty}}, \qquad\qquad
\|e_{ u^{(i)}}\|_{\infty}=\frac{\|u^{(i)}- u_{h}^{(i)}\|_{\infty}}{\|  u^{(i)}\|_{\infty}},\quad
\quad i = 1,2,3,$$
and 
$$\|e_{p}\|_2=\frac{\|p-p_h\|_2}{\|p\|_2},\qquad\qquad
\|e_{ u^{(i)}}\|_2=\frac{\|u^{(i)}-u_{h}^{(i)}\|_2}{\|u^{(i)}\|_2},\qquad\;
\quad i = 1,2,3.$$
In all the examples,  the GMRES iterative method is employed to solve the discrete boundary integral equations. 
The GMRES iteration starts with a zero initial guess and stops when the iterated residual in the discrete $\ell^2$-norm relative to that of the initial residual is less than a prescribed tolerance $\epsilon=10^{-9}$. Furthermore, the corresponding bounding cube for the interface problem is set to $\mathcal{B}=[-1.2,1.2]^3$.

\subsection{Examples for Stokes problems}
\begin{figure}[!ht]
  \centering
  \subfigure[$u_h^{(1)}$]{
    \includegraphics[width=1.35in]{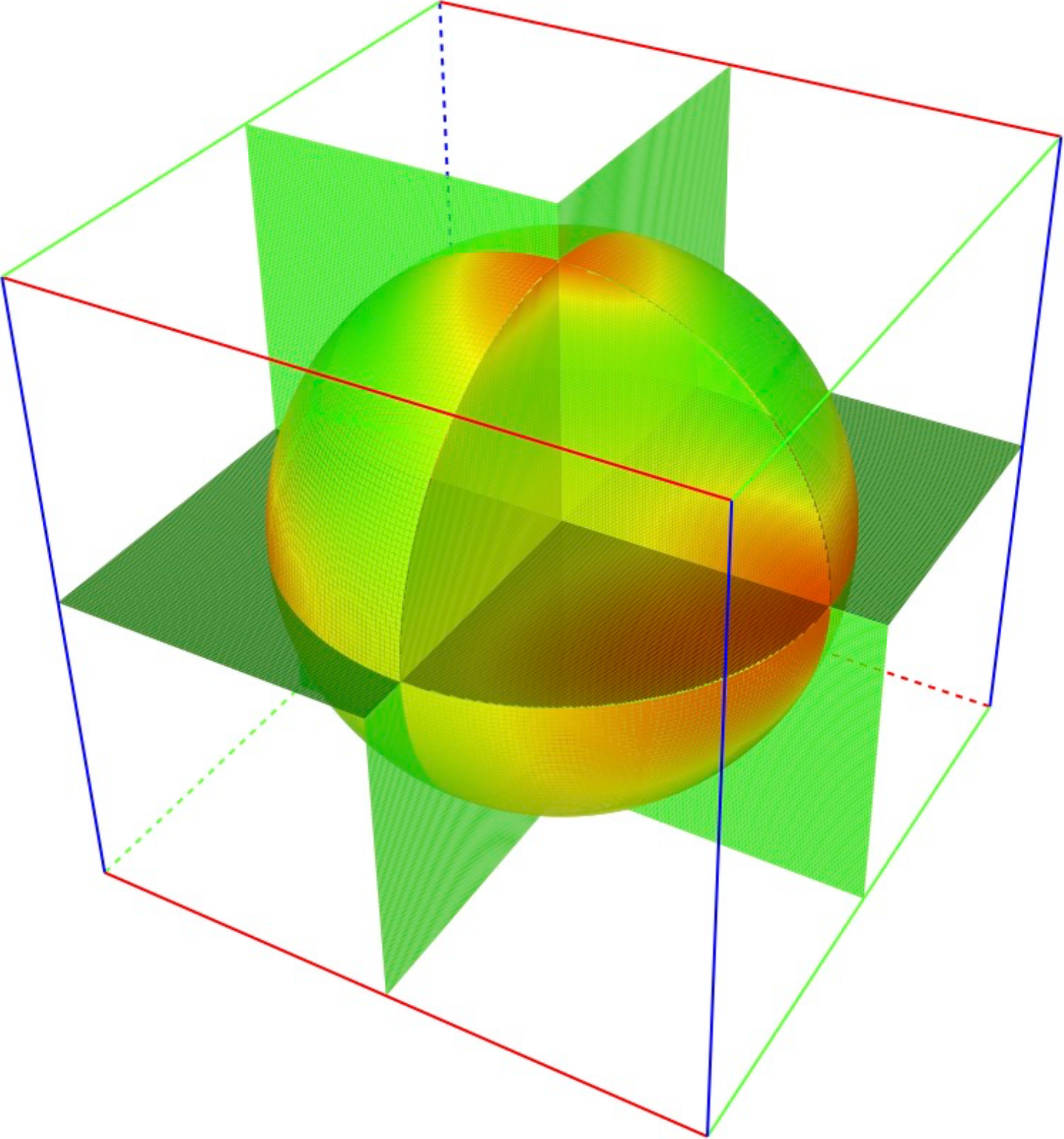}
  }
  \subfigure[$u_h^{(2)}$]{
    \includegraphics[width=1.35in]{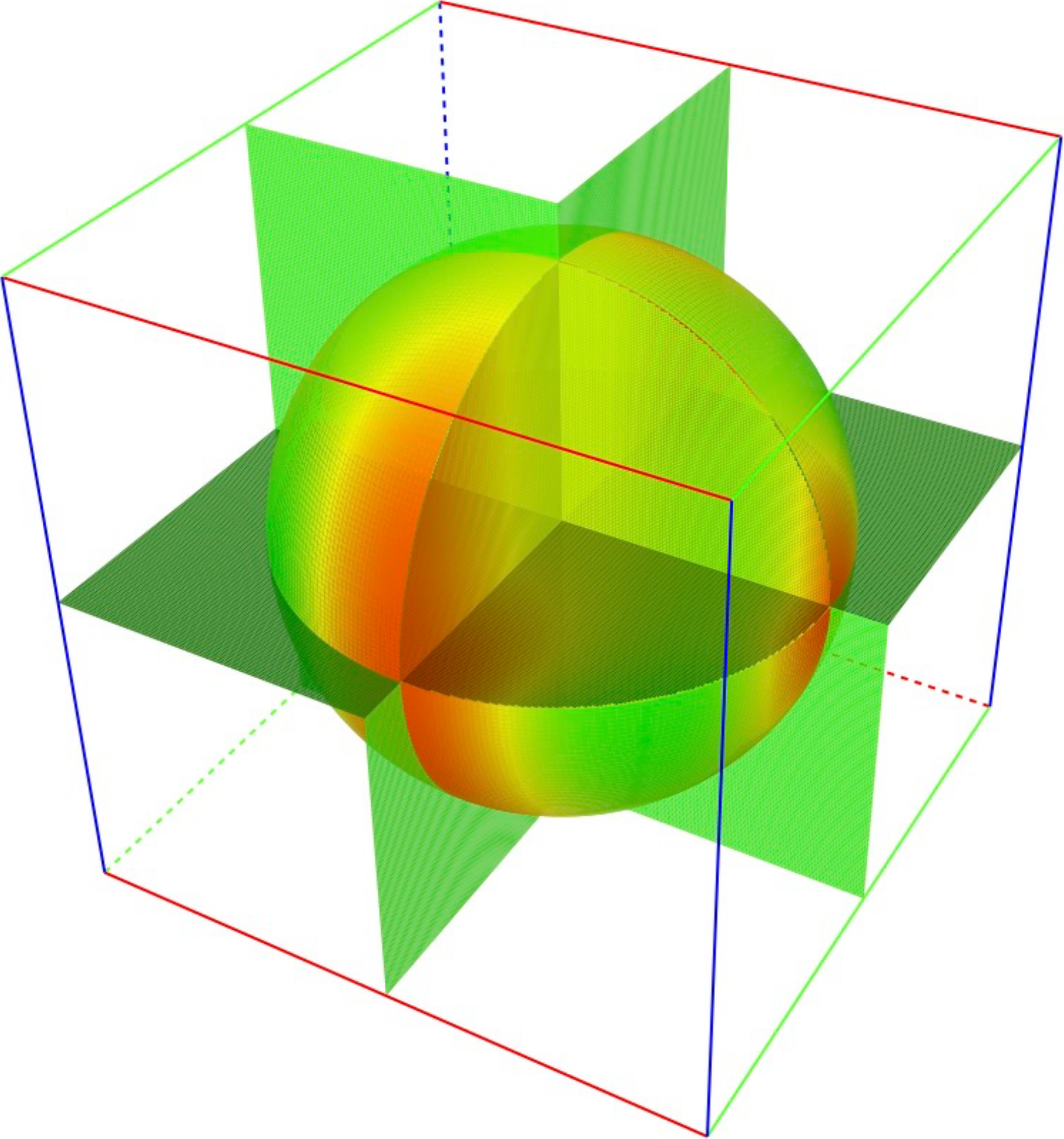}
  }
  \subfigure[$u_h^{(3)}$]{
    \includegraphics[width=1.35in]{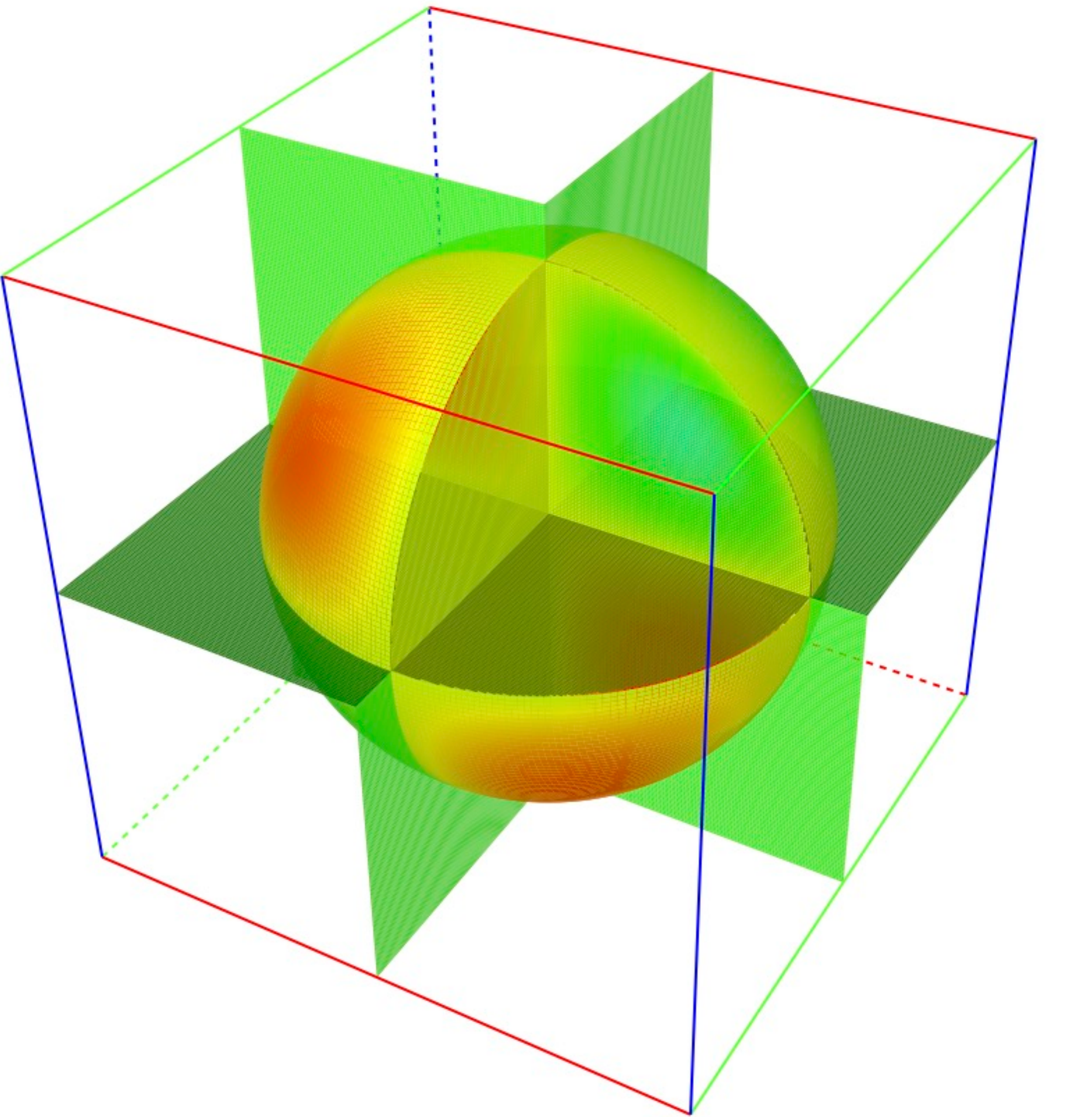}
  }
  \subfigure[$p_h$]{
    \includegraphics[width=1.35in]{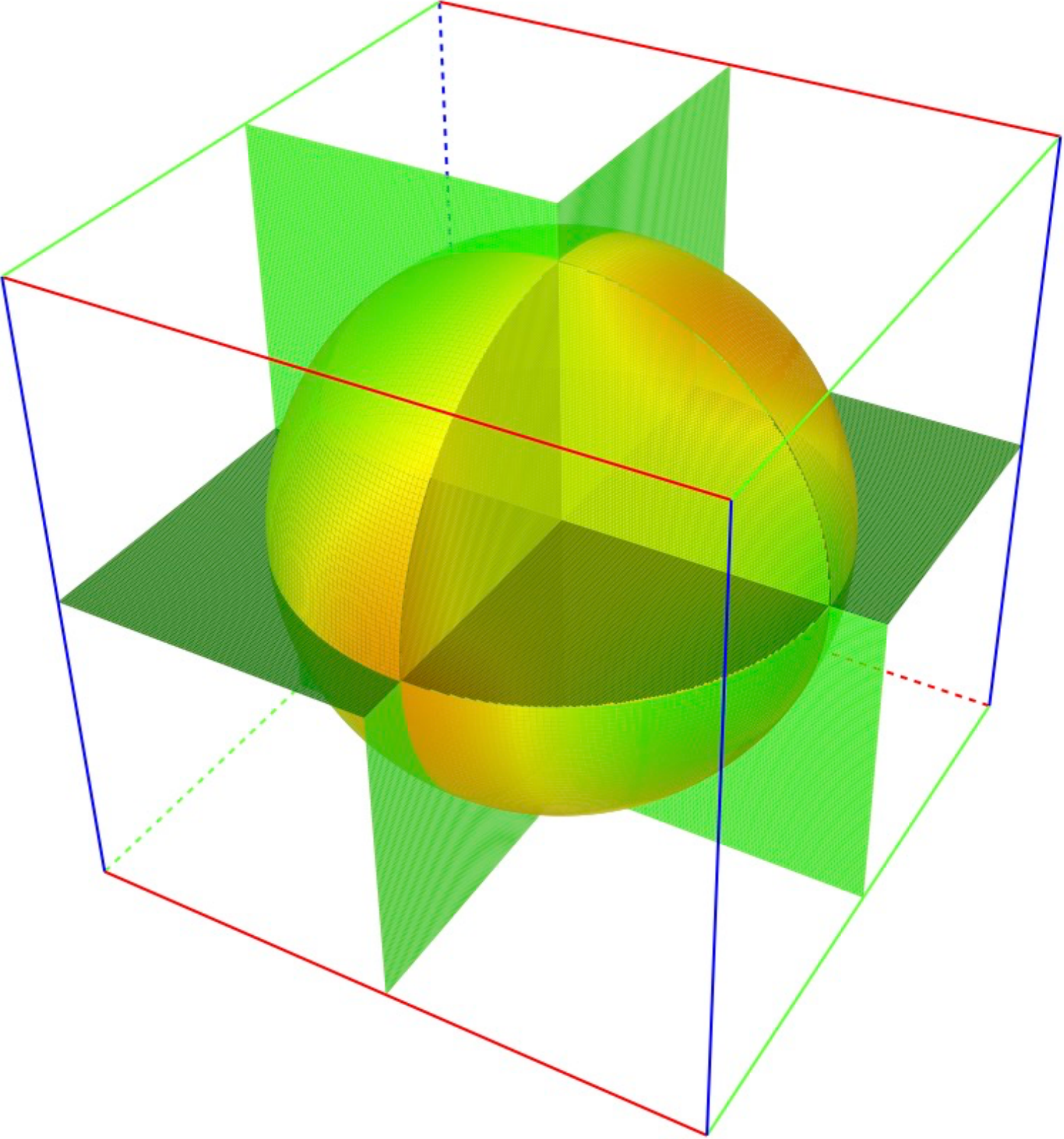}
  }
  \caption{ Numerical solution $\vect u_h$ and $p_h$ in a sphere for Stokes problem.}
\label{figEX1}
\end{figure}

\begin{figure}[!ht]
  \centering
  \subfigure[$u_h^{(1)}$]{
    \includegraphics[width=1.35in]{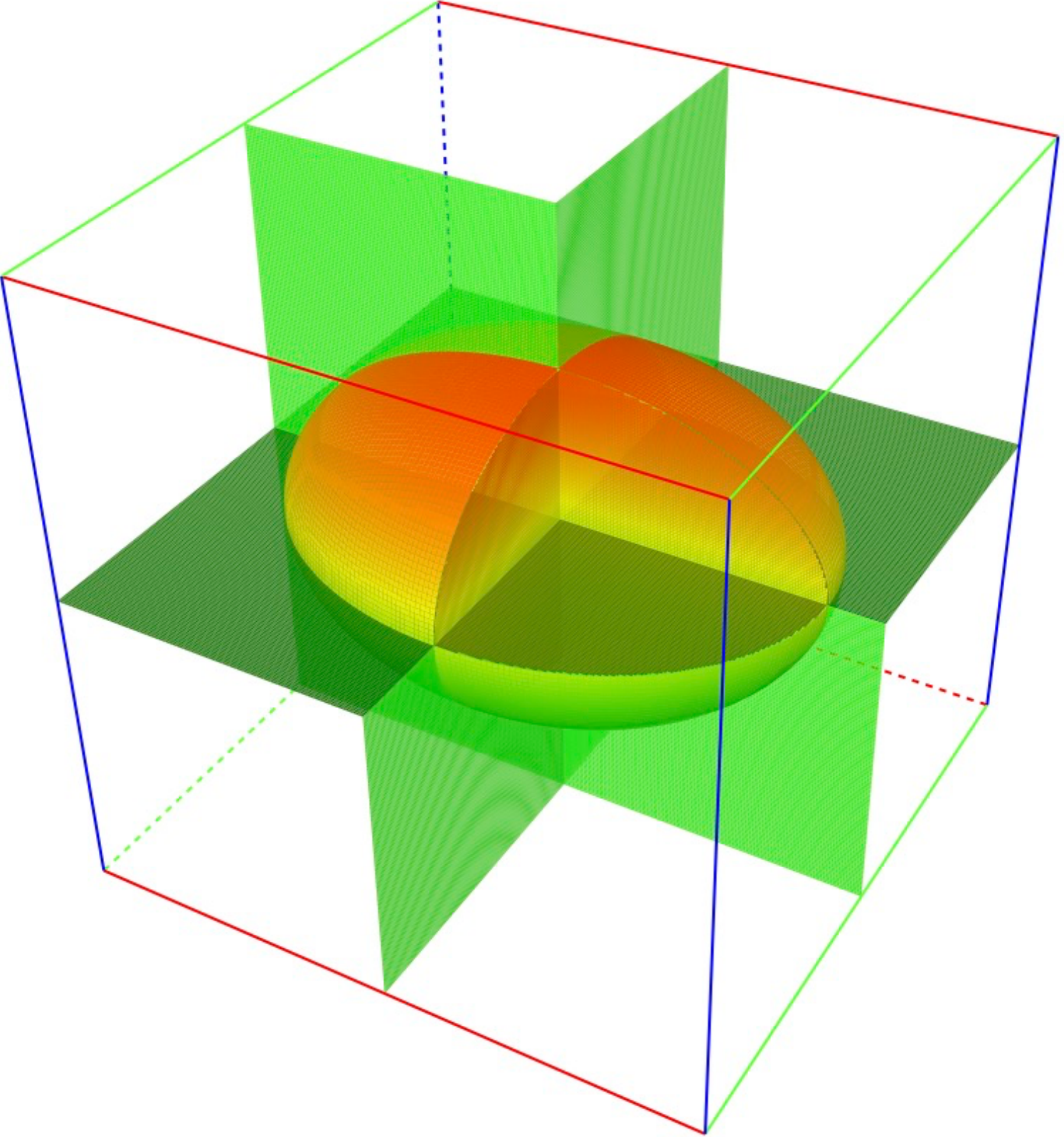}
  }
  \subfigure[$u_h^{(2)}$]{
    \includegraphics[width=1.35in]{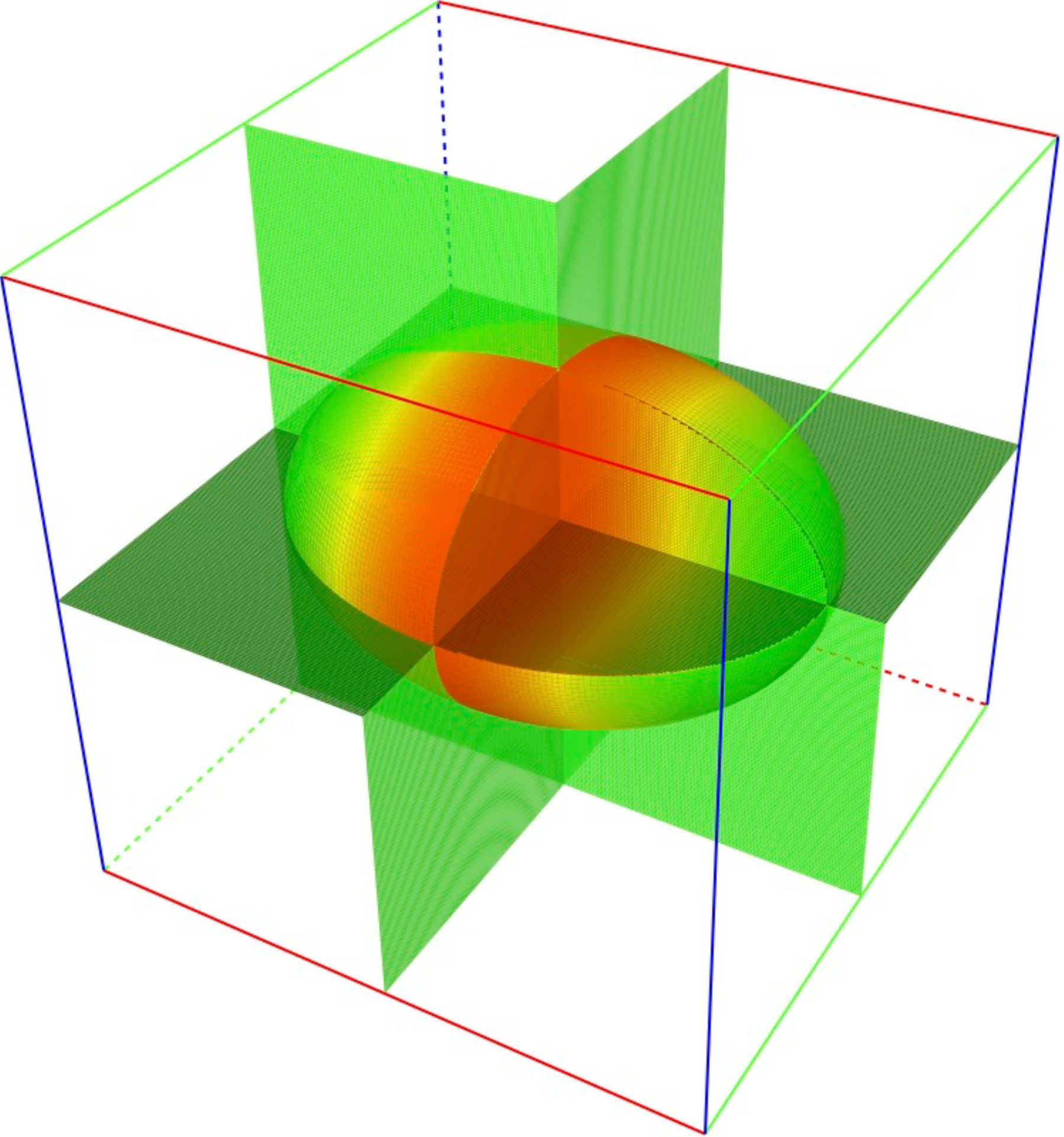}
  }
  \subfigure[$u_h^{(3)}$]{
    \includegraphics[width=1.35in]{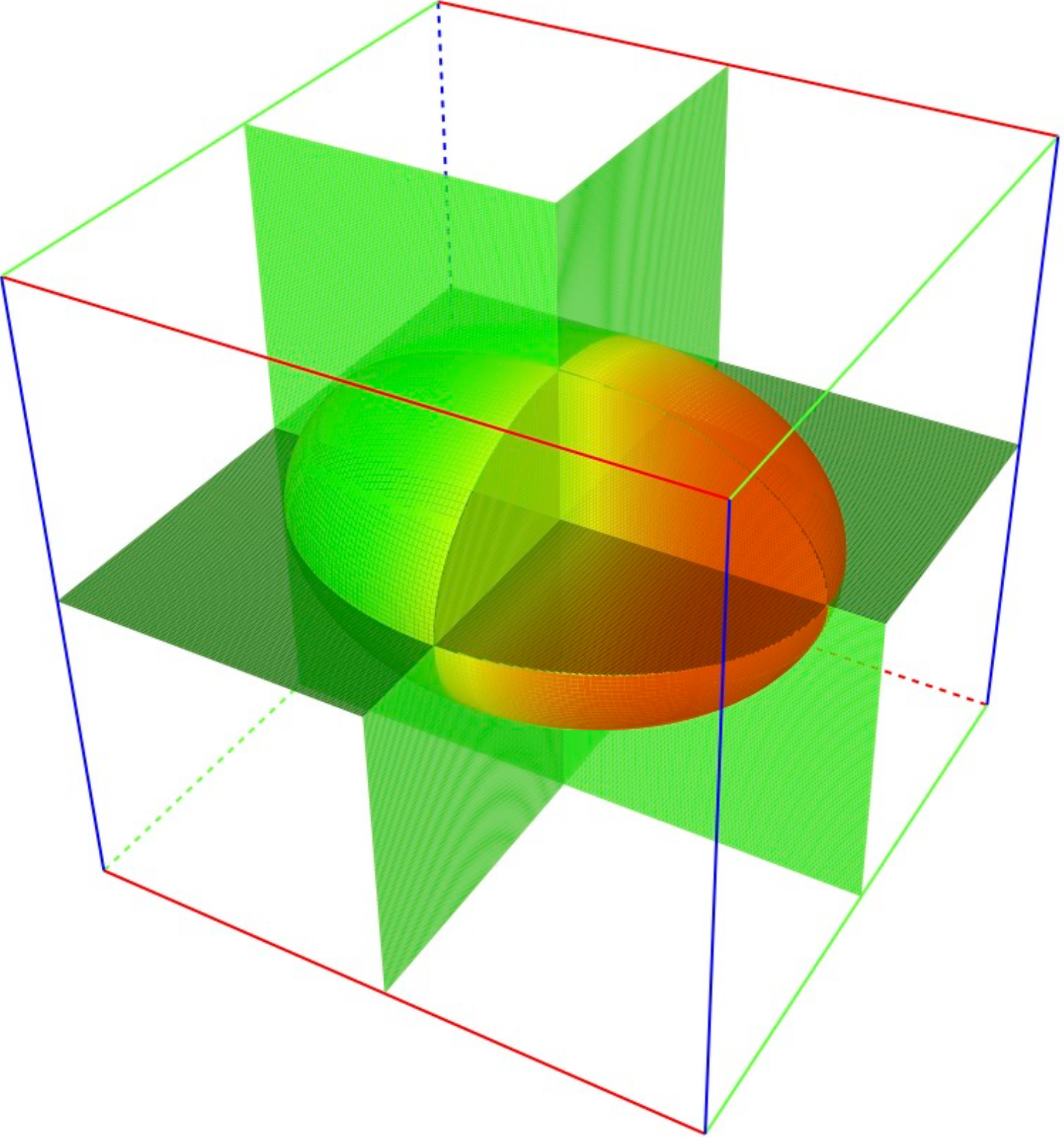}
  }
  \subfigure[$p_h$]{
    \includegraphics[width=1.35in]{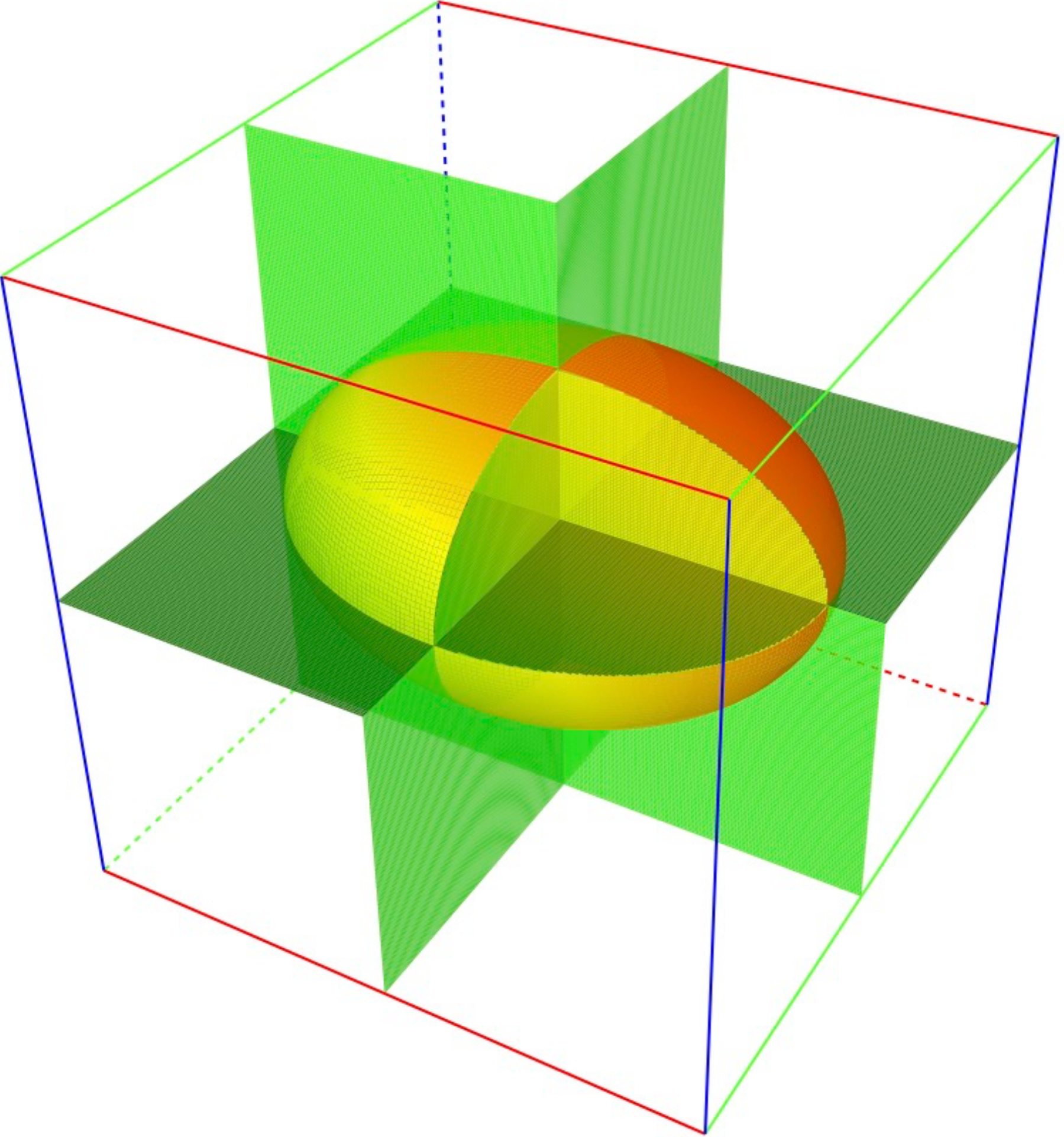}
  }
  \caption{ Numerical solution $\vect u_h$ and $p_h$ in an ellipsoid for Stokes problem.}
\label{figEX2}
\end{figure}

\begin{figure}[!ht]
  \centering
  \subfigure[$u_h^{(1)}$]{
    \includegraphics[width=1.35in]{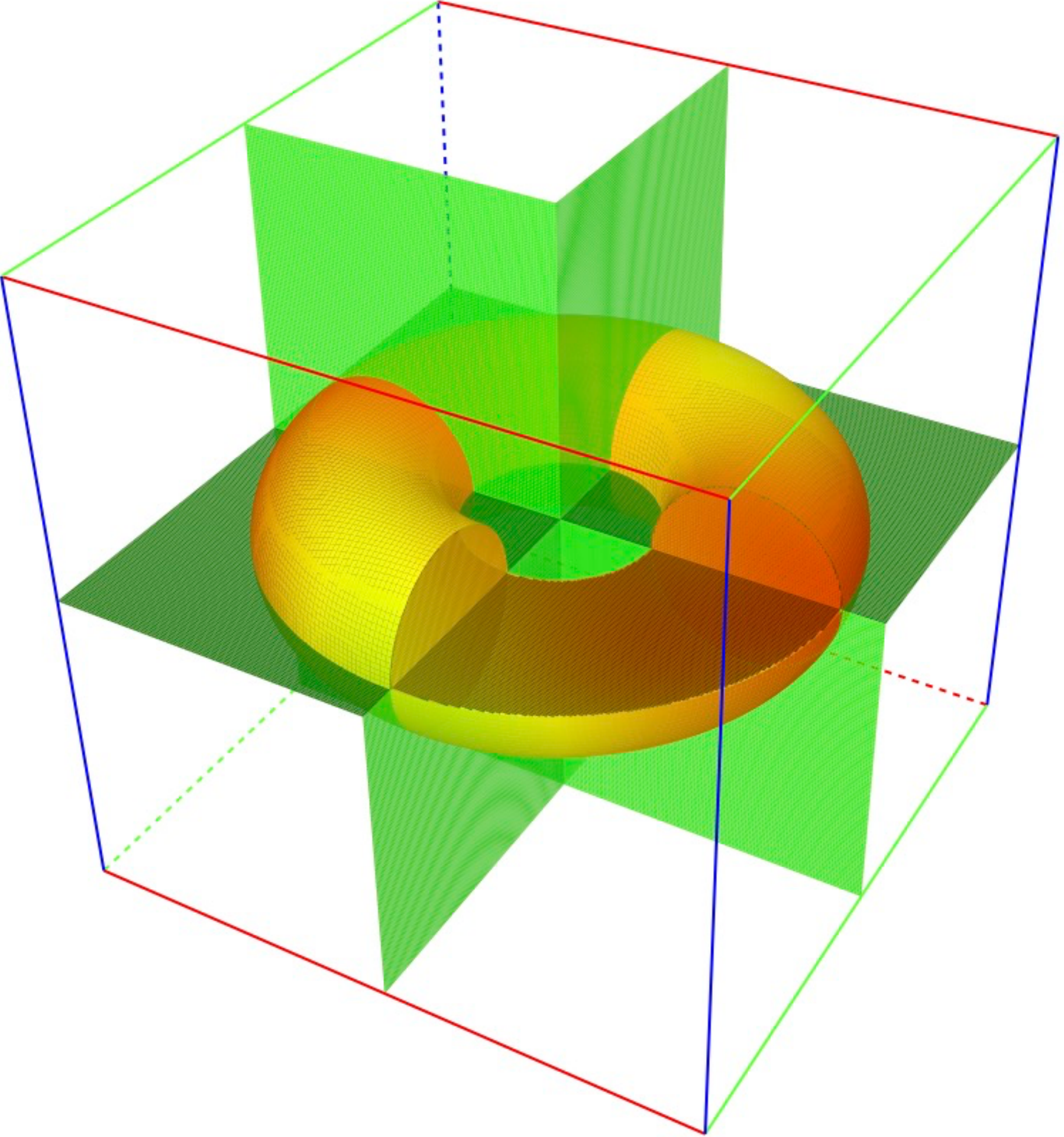}
  }
  \subfigure[$u_h^{(2)}$]{
    \includegraphics[width=1.35in]{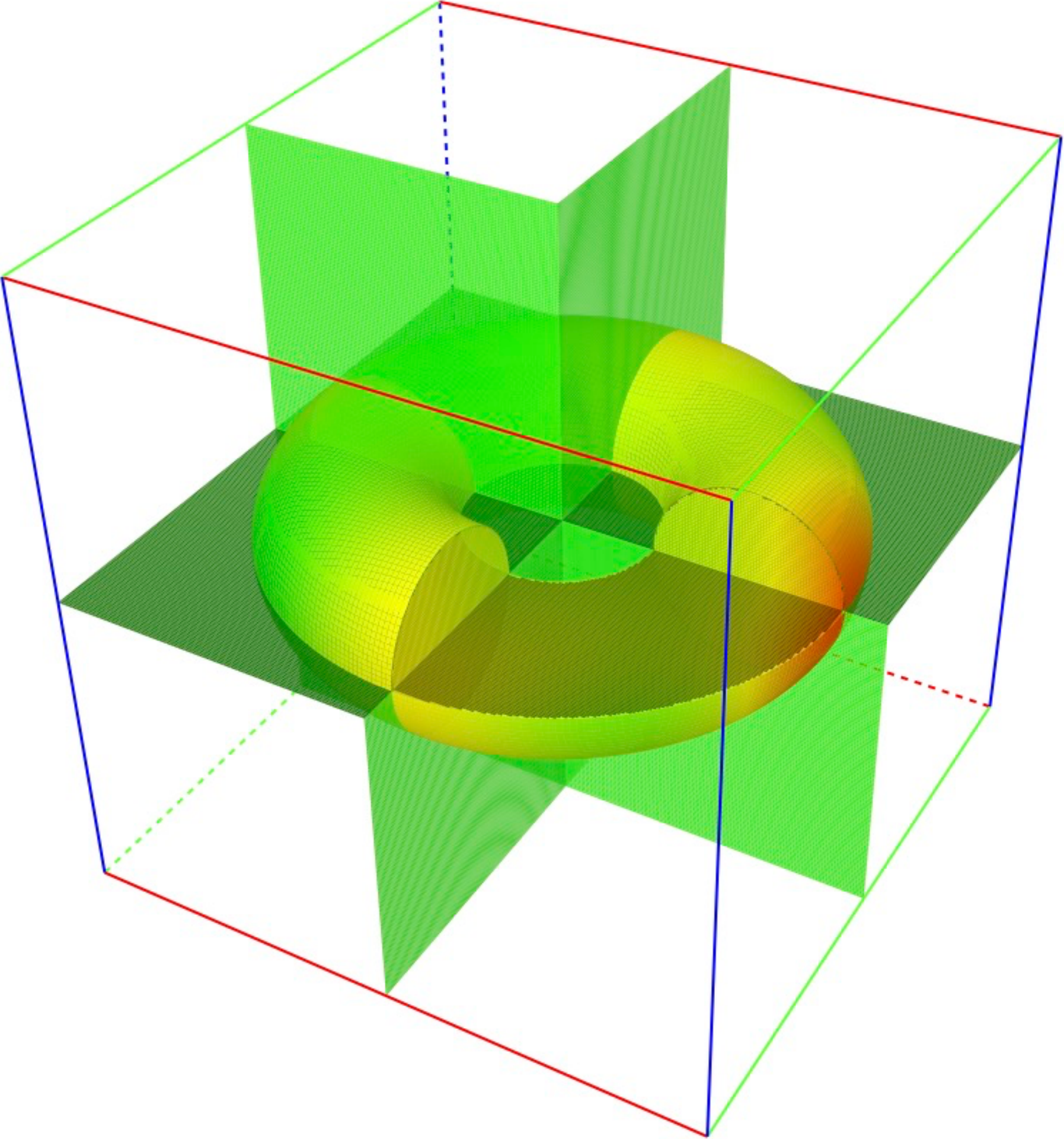}
  }
  \subfigure[$u_h^{(3)}$]{
    \includegraphics[width=1.35in]{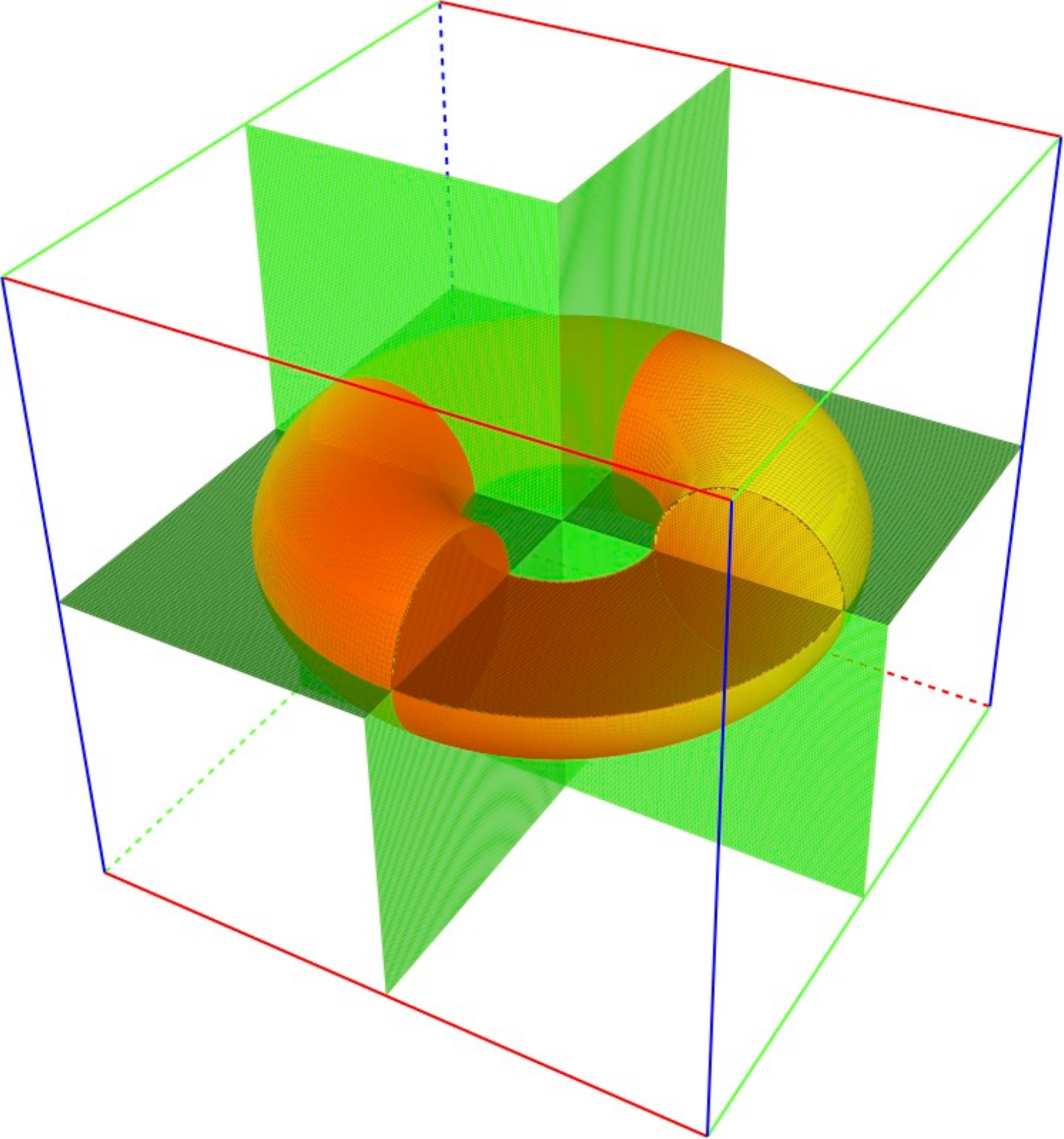}
  }
  \subfigure[$p_h$]{
    \includegraphics[width=1.35in]{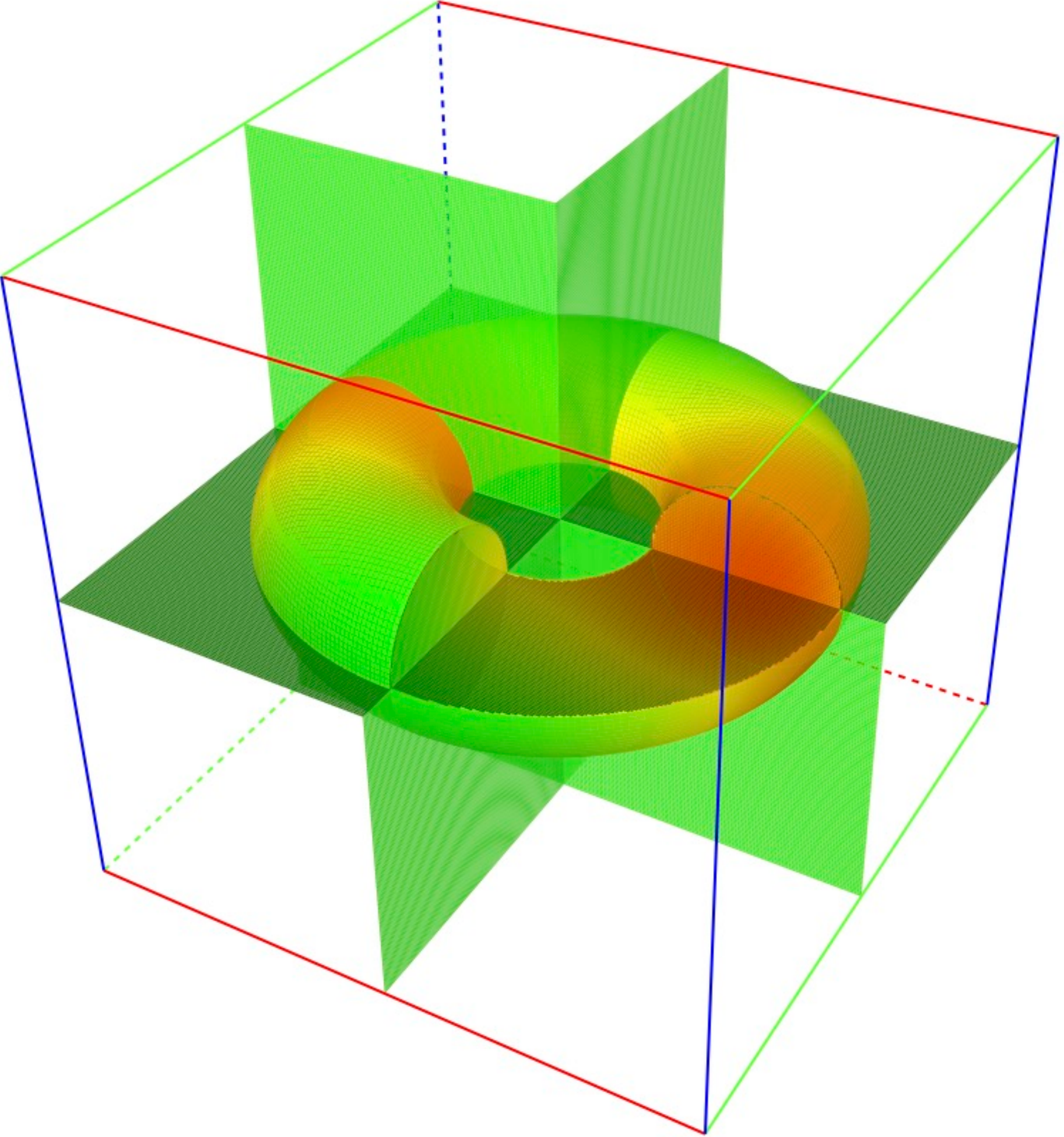}
  }
  \caption{ Numerical solution $\vect u_h$ and $p_h$ in a torus for Stokes problem.}
\label{figEX3}
\end{figure}

Three Stokes problems on different 3D irregular domains are considered in this subsection. Numerical results are listed in Tables \ref{sphere1M}-\ref{tabtorus3L}. The grid sizes are listed in the first column, and the number of GMRES iterations is recorded in the second column. The third to the last columns show the errors of the numerical solution in discrete maximum and $\ell^2$-norm, as well as the convergence rates. The color mapped numerical solution $\vect u_h$ and $p_h$ on a $256\times 256$ mesh are shown in Fig. \ref{figEX1}-\ref{figEX3}. It is noted that a fine mesh is used in all plots presented in this section for a better resolution of the irregular boundary.


{\bf Example 5.1.} This example solves the Stokes Dirichlet problem on a sphere with radius $r = 1$, which is located at the origin of the coordinates. The exact solution is given by 
\begin{equation*}
\begin{cases}
u^{(1)} =x^2(2-x^2)+4xy(x^2+y^2-1)+z^2(3z^2-6x^2-2),\\[4pt]
u^{(2)} = x^2(3x^2-6y^2-1)+y^2(1-y^2),\\[4pt]
u^{(3)} =4zx(z^2+x^2-1),\\[4pt]
\;\quad p = 8y(3x^2-y^2)+8x(3z^2-x^2).
\end{cases}
\end{equation*}
The external forcing term $\vect f$ and the boundary condition $\vect g$ can be determined from the exact solution. The errors and convergence rates in the maximum norms and $\ell^2$ norms are shown in Table \ref{sphere1M} and Table \ref{sphere1L} respectively, which indicate that the velocity is of second-order accuracy in both the discrete maximum and $\ell^2$ norm, and the pressure is also second-order accurate in the $\ell^2$ norm. The numbers of GMRES iterations are also presented, which shows that the number of GMRES iterations is independent of the size of the mesh. 

\begin{table}[!ht]
\caption{Maximum error and convergence rates of {\em Example} 1  }
\begin{center}
\vspace{-0.1cm}
\begin{tabular}{|c|c|c|c|c|c|c|c|c|c|}
\hline
N & \# step  & $\|e_{\vect u_1}\|_{\infty}$
& rate & $\|e_{\vect u_2}\|_{\infty}$ & rate &$\|e_{\vect u_3}\|_{\infty}$ & rate & $\|e_{p}\|_{\infty}$ & rate \\
\hline
64 & 12 & 4.80e-3 & - & 3.65e-3 & - & 9.23e-3 & - & 2.33e-2 & -\\
128 & 11 & 5.59e-4 & 3.10 & 6.84e-4 & 2.42 & 1.39e-3 & 2.73 & 6.40e-3 & 1.86\\
256 & 11 & 9.55e-5 & 2.55 & 8.67e-5 & 2.98 & 1.31e-4 & 3.41 & 1.91e-3 & 1.74\\
512 & 11 & 1.34e-5 & 2.83 & 1.44e-5 & 2.59 & 1.78e-5 & 2.88 & 6.06e-4 & 1.66\\
\hline
\end{tabular}
\end{center}
\label{sphere1M}
\end{table}

\begin{table}[!ht]
\caption{$\ell^2$-error and convergence rates of {\em Example} 1  }
\begin{center}
\vspace{-0.1cm}
\begin{tabular}{|c|c|c|c|c|c|c|c|c|}
\hline
grid size &  $\|e_{\vect u_1}\|_{2}$
& rate & $\|e_{\vect u_2}\|_{2}$ & rate &$\|e_{\vect u_3}\|_{2}$ & rate & $\|e_{p}\|_{2}$ & rate \\
\hline
64 & 1.06e-3 & - & 9.16e-4 & - & 1.40e-3 & - & 5.46e-3 & -\\
128 &  9.12e-5 & 3.54 & 8.73e-5 & 3.39 & 1.18e-4 & 3.57 & 9.56e-4 & 2.51 \\
256 & 1.12e-5 & 3.03 & 9.72e-6 & 3.17 & 1.57e-5 & 2.91 & 2.04e-4 & 2.23\\
512 &  1.43e-6 & 2.97 & 1.37e-6 & 2.83 & 1.92e-6 & 3.03 & 3.82e-5 & 2.42\\
\hline
\end{tabular}
\end{center}
\label{sphere1L}
\end{table}

{\bf Example 5.2.} In this example, the computational domain is an ellipsoid $$\Big\{(x,y,z)\ |\ x^2+\dfrac{y^2}{0.64}+\dfrac{z^2}{0.36} = 1.\Big\}$$ Consider the exact solution  
\begin{equation*}
\begin{cases}
u^{(1)} =2x^3y/3-z^2x^2+e^z+\sin(\pi z),\\[4pt]
u^{(2)} =-x^2y^2+\cos(\pi x),\\[4pt]
u^{(3)} =2z^3x/3+e^x+\sin(\pi x),\\[4pt]
\;\quad p = (x-1)^3(y-1)^3(z-1)^3.
\end{cases}
\end{equation*}
Normalized errors for the velocity components $u^{(1)}$, $u^{(2)}$, $u^{(3)}$ and the pressure $p$ in the discrete maximum norm and the $\ell^2$ norm are listed separately in Table \ref{ellipse2M} and \ref{ellipse2L}. It can be seen that the components of velocity are all second-order accurate in both the discrete maximum norm and discrete  $\ell^2$  norm, and the pressure is second accurate in $\ell^2$ norm. The GMRES iteration number for different mesh sizes is also shown in Table \ref{ellipse2M}-\ref{ellipse2L}.
\begin{table}[!ht]
\caption{Maximum error and convergence rates of {\em Example} 2  }
\begin{center}
\vspace{-0.1cm}
\begin{tabular}{|c|c|c|c|c|c|c|c|c|c|}
\hline
N & \# step  & $\|e_{\vect u_1}\|_{\infty}$
& rate & $\|e_{\vect u_2}\|_{\infty}$ & rate &$\|e_{\vect u_3}\|_{\infty}$ & rate & $\|e_{p}\|_{\infty}$ & rate \\
\hline
128 & 14 & 5.44e-4 & - & 7.65e-4 & - & 2.53e-4 & - & 2.22e-2 & - \\
256 & 15 & 5.94e-5 & 3.20 & 8.01e-5 & 3.26 & 4.13e-5 & 2.61 & 7.05e-3 & 1.65 \\
512 & 14 & 8.44e-6 & 2.82 & 2.07e-5 & 1.95 & 6.80e-6 & 2.60 & 1.77e-3 & 1.99\\
\hline
\end{tabular}
\end{center}
\label{ellipse2M}
\end{table}

\begin{table}[!ht]
\caption{$\ell^2$-error and convergence rates of {\em Example} 2  }
\begin{center}
\vspace{-0.1cm}
\begin{tabular}{|c|c|c|c|c|c|c|c|c|}
\hline
N &  $\|e_{\vect u_1}\|_{2}$
& rate & $\|e_{\vect u_2}\|_{2}$ & rate &$\|e_{\vect u_3}\|_{2}$ & rate & $\|e_{p}\|_{2}$ & rate \\
\hline
128 &  2.95e-5 & - & 6.43e-5 & - & 1.82e-5 & - & 1.14e-3 & -\\
256 &  6.02e-6 & 2.29 & 1.26e-5 & 2.35 & 4.18e-6 & 2.12 & 2.94e-4 & 1.96 \\
512 &  1.29e-6 & 2.22 & 2.07e-6 & 2.61 & 9.42e-7 & 2.15 & 9.10e-5 & 1.69\\
\hline
\end{tabular}
\end{center}
\label{ellipse2L}
\end{table}

{\bf Example 5.3.} In this example, a more complicated domain is considered, which is given by $$\Omega = \{(x, y, z)\in \mathcal{R}^3 : (c-\sqrt{(x^2+y^2)})^2+z^2<a^2\}$$ 
with $a = 0.35$ and $c = 0.7$.  The exact solution is determined by:
\begin{equation*}
\begin{cases}
u^{(1)} = -4xy(1-x^2-y^2) -x^2(x^2+6z^2-2)+z^2(3z^2-2)\\[4pt]
\quad\quad + \exp(\cos y) + \exp(\sin z) ,\\[4pt]
u^{(2)} = x^2(3x^2-6y^2-2)-y^2(y^2-2)+\exp(\sin x),\\[4pt]
u^{(3)} = -4(1-x^2-z^2)xz + \exp(\cos x),\\[4pt]
\quad\;p = \exp(1-y^2-z^3)\sin(x^2+1).
\end{cases}
\end{equation*}
The maximum errors, $\ell^2$-errors and the corresponding convergence rates are displayed in Table \ref{tabtorus3M}-\ref{tabtorus3L}, which indicate that second-order accuracy of the solutions $(\vect u_h, p_h)$ is achieved. And the number of GMRES iterations shown in the second column of the tables is almost independent of the mesh size.
 \begin{table}[!ht]
\caption{Maximum error and convergence rates of {\em Example} 3  }
\begin{center}
\vspace{-0.1cm}
\begin{tabular}{|c|c|c|c|c|c|c|c|c|c|}
\hline
N & \# step  & $\|e_{\vect u_1}\|_{\infty}$
& rate & $\|e_{\vect u_2}\|_{\infty}$ & rate &$\|e_{\vect u_3}\|_{\infty}$ & rate  & $\|e_{p}\|_{\infty}$ & rate \\
\hline
128 & 24 & 8.97e-4 & - &  8.29e-4 & - & 7.42e-4 & - & 5.54e-1 & -\\
256 & 26 & 7.54e-5 & 3.57 & 5.80e-5 & 3.84 & 6.19e-5 & 3.58 & 1.00e-1 & 2.47\\
512 & 26 & 1.01e-5 & 2.90 & 9.86e-6 & 2.56 & 1.32e-5 & 2.23 & 1.86e-2 & 2.43\\
\hline
\end{tabular}
\end{center}
\label{tabtorus3M}
\end{table}

\begin{table}[!ht]
\caption{$\ell^2$-error and  its  convergence rates of {\em Example} 3  }
\begin{center}
\vspace{-0.1cm}
\begin{tabular}{|c|c|c|c|c|c|c|c|c|}
\hline
N &  $\|e_{\vect u_1}\|_{2}$
& rate & $\|e_{\vect u_2}\|_{2}$ & rate &$\|e_{\vect u_3}\|_{2}$ & rate &  $\|e_{p}\|_{2}$ & rate \\
\hline
128 & 2.22e-5 & - & 7.64e-5 & - & 2.62e-5 & - & 1.50e-2 & - \\
256 &  4.74e-6 & 2.23 & 1.26e-5 & 2.60 & 3.95e-6 & 2.73 & 2.94e-3 & 2.35 \\
512 &  1.03e-6 & 2.20 & 2.50e-6 & 2.33 & 9.14e-7 & 2.11 & 8.82e-4 & 1.74 \\
\hline
\end{tabular}
\end{center}
\label{tabtorus3L}
\end{table}

\subsection{Examples for Navier problems}
To test the efficiency of the KFBI for the Navier boundary value problems, three examples are considered here.
Numerical results are listed in Tables \ref{tabsphere-max}-\ref{tabtorus-l2}. Each table has eight columns, showing the size of the Cartesian grid, the numbers of the GMRES iteration in solving BIE,  the errors of the numerical solution in the maximum norm or discrete $\ell^2$-norm, as well as the convergence rates.  Fig. 4-6 show the color mapped numerical solution $\vect u_h$ of each example on a $256\times 256$ mesh respectively. 

{\bf Example 5.4.} 
In this example, the computational domain is a sphere, which is given by $x^2+y^2+z^2=1$.  The exact solution is determined by:
\begin{equation*}
\begin{cases}
u^{(1)} = \sin(2\pi r^2)\cos(2y)\cos(z),\\[3pt]
u^{(2)}= \sin(2\pi r^2)\cos(x)\cos(x+z),\\[3pt]
u^{(3)} = \sin(2\pi r^2)\cos(\pi x)\cos(y).
\end{cases}
\end{equation*}
The material data are chosen as $E=1000$, and $\nu=1/10$. The maximum errors and the corresponding convergence rates are displayed in Table \ref{tabsphere-max}.
The $\ell^2$-errors and the corresponding convergence rates are shown in Table \ref{tabsphere-l2}. As expected, the proposed method is second-order accurate for the numerical solutions. Additionally, the number of iterations does not depend on the mesh size and it is a relatively small number. 

\begin{table}[!ht]
\caption{Maximum-error and  its  convergence rates of Example 4.4 }
\begin{center}
\vspace{-0.1cm}
\begin{tabular}{|c|c|c|c|c|c|c|c|}
\hline
N & \# step  & $\|e_{\vect u_1}\|_{2}$
& rate & $\|e_{\vect u_2}\|_{2}$ & rate &$\|e_{\vect u_3}\|_{2}$ & rate\\
\hline
64   &  14  &1.50e-2  & -        & 1.33e-2  & -        & 1.60e-2 &      -  \\
128   &  14  &4.09e-3  &  1.87 & 3.03e-3  & 2.13  & 4.02e-3 & 1.99 \\
256   &  14  &9.97e-4  &  2.04 & 7.98e-4  & 1.92  & 9.81e-4 & 2.03 \\
512 &  15  &2.62e-4  &  1.92 & 2.10e-4  & 1.93	& 2.60e-4 & 1.92\\	
\hline
\end{tabular}
\end{center}
\label{tabsphere-max}
\end{table}

\begin{table}[!ht]
\caption{$\ell^2$-error and  its  convergence rates of Example 4.4 }
\begin{center}
\vspace{-0.1cm}
\begin{tabular}{|c|c|c|c|c|c|c|}
\hline
N & $\|e_{\vect u_1}\|_{2}$
& rate & $\|e_{\vect u_2}\|_{2}$ & rate &$\|e_{\vect u_3}\|_{2}$ & rate\\
\hline
64   &  9.03e-3  & -        & 8.03e-3  & -        & 9.63e-3 &      -  \\
128   &  2.12e-3  &  2.09 & 1.90e-3  & 2.08  & 2.29e-3 & 2.07 \\
256   &  4.05e-4  &  2.39 & 3.84e-4  & 2.31  & 4.97e-4 & 2.20 \\
512 & 1.01e-4  &  2.00 & 9.65e-5  & 1.99	& 1.25e-4 & 1.99\\	
\hline
\end{tabular}
\end{center}
\label{tabsphere-l2}
\end{table}

\begin{figure}[!ht]
  \centering
  \subfigure[$u_h^{(1)}$]{
    \includegraphics[width=1.5in]{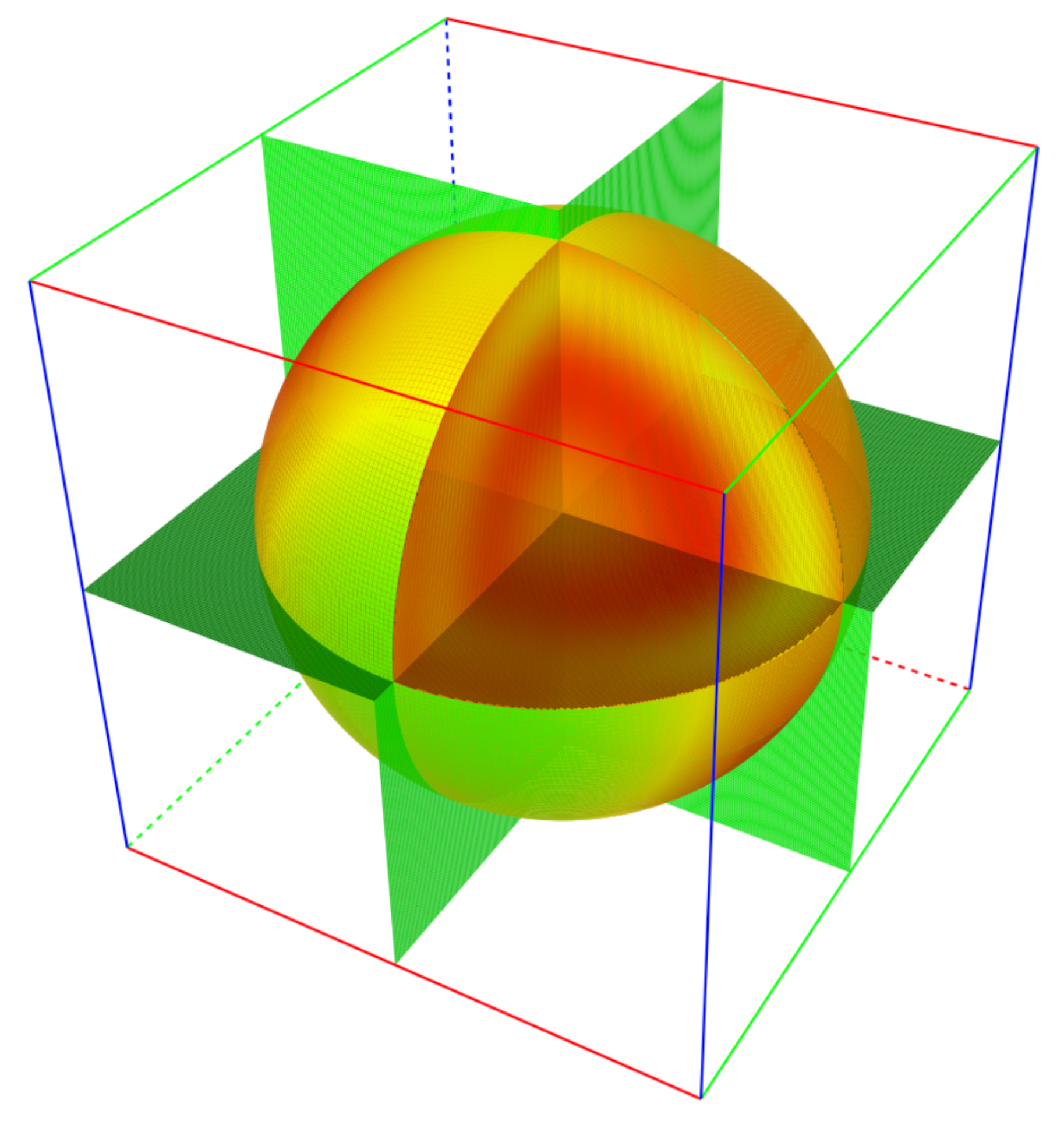}
  }
  \subfigure[$u_h^{(2)}$]{
    \includegraphics[width=1.5in]{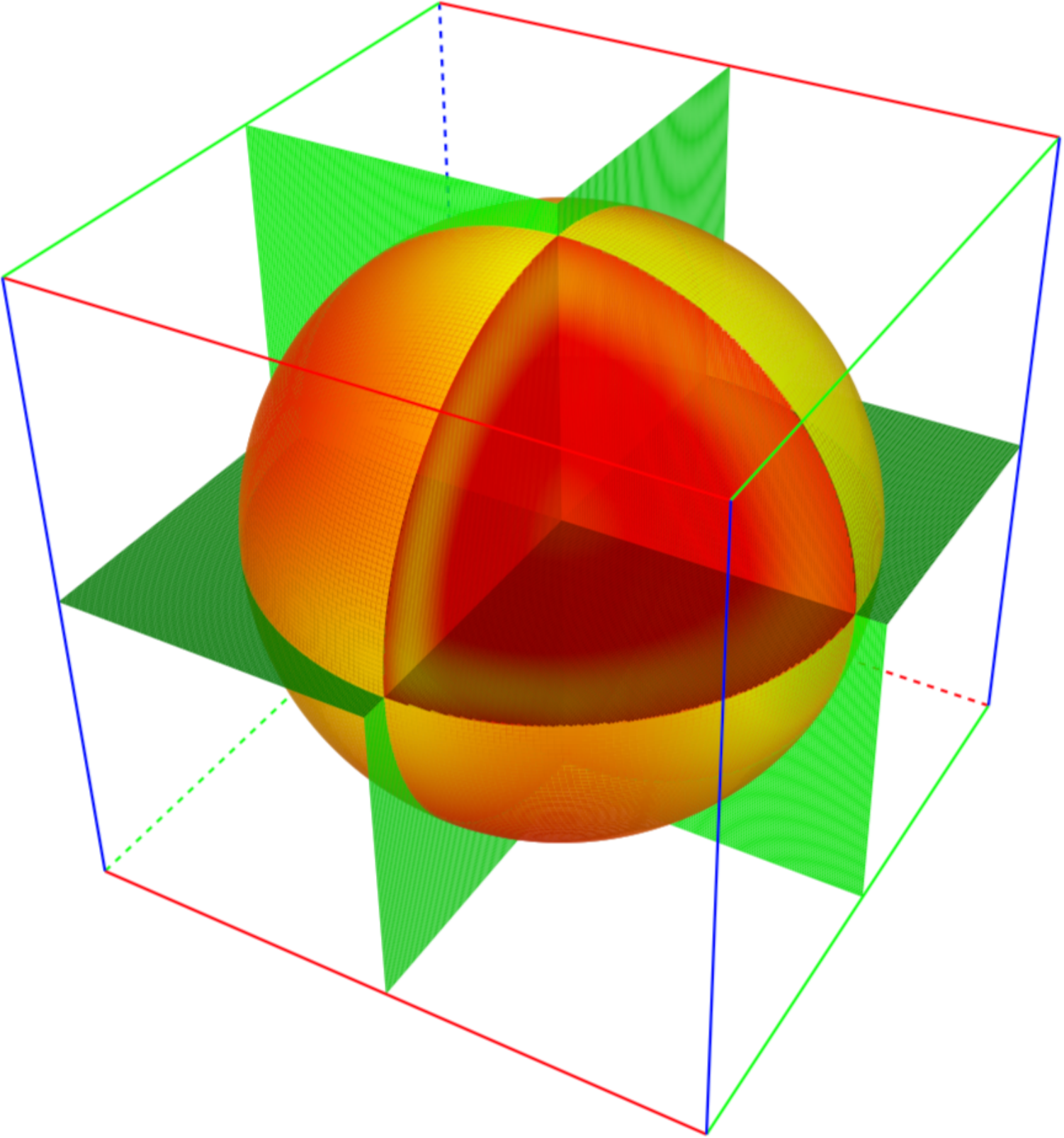}
  }
    \subfigure[$u^{(3)}_h$]{
    \includegraphics[width=1.5in]{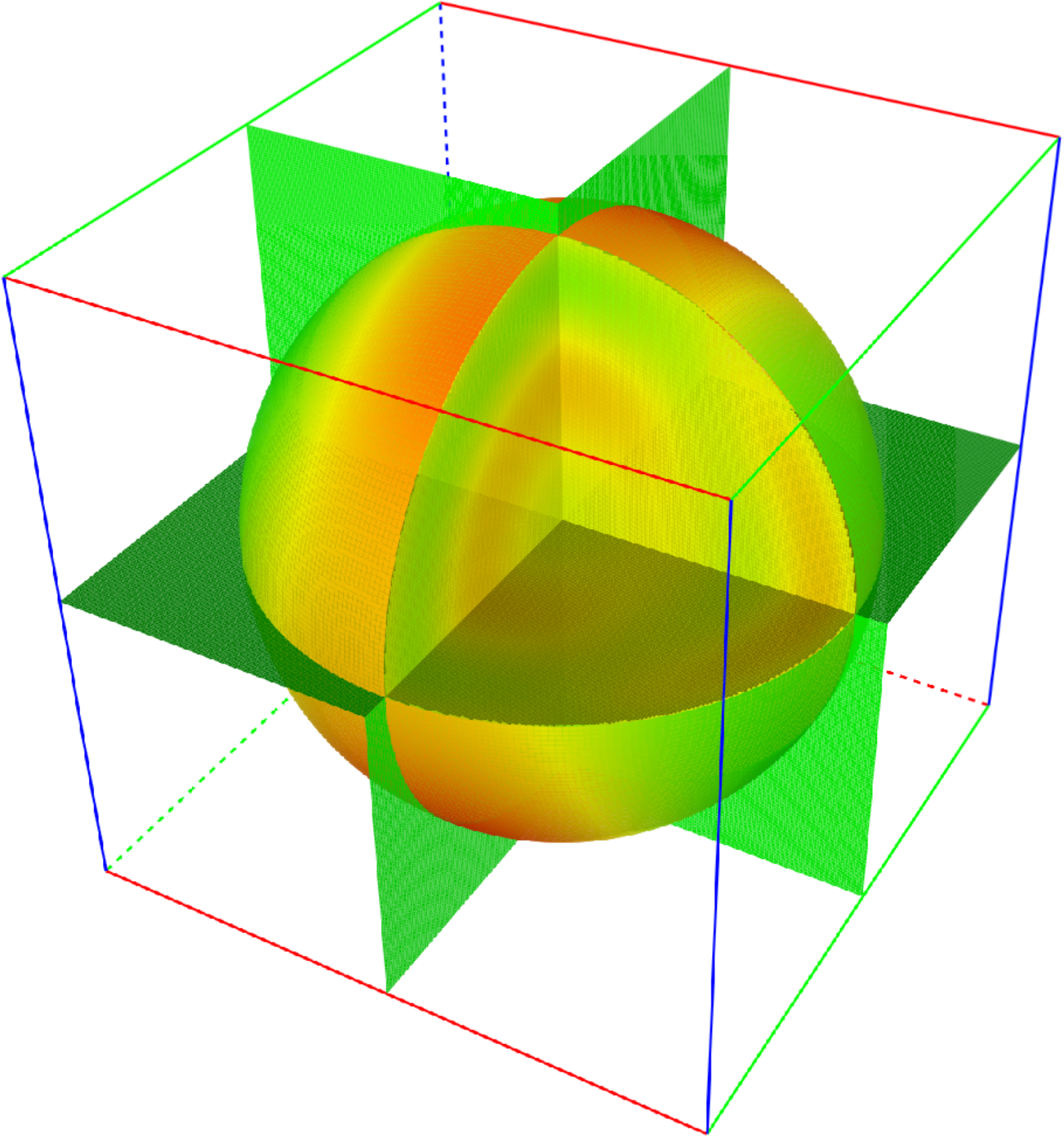}
  }
\caption{ Numerical solution $\vect u_h$ in a sphere for Navier problem.}
\label{figsphere}
\end{figure}

{\bf Example 5.5.}
This example solves the Navier equations on an ellipsoid domain, which is determined by $x^2/r_a^2+y^2/r_b^2+z^2/r_c^2=1$ with $r_a=1, r_b=0.8, r_c=0.6$. The exact solution is given by
\begin{equation*}
\begin{cases}
u^{(1)}=x^2+y^2+z^2-4+\cos x\cos y\cos z,\\[4pt]
u^{(2)}=x^2+y^2+z^2-4+xy+\cos x\cos y\cos z, \\[4pt]
u^{(3)}=x^2+y^2+z^2-4+yz+\cos x\cos y\cos z.
\end{cases}
\end{equation*}
The material data are chosen as $E=4000$, and $\nu=1/20$. 
Numerical results with second-order accuracy are summarized in Tables \ref{tabellipse-max} and \ref{tabellipse-l2}. The efficiency and accuracy of the numerical solution observed here are consistent with Example 4.4. 

\begin{table}[!ht]
\caption{Maximum-error and  its  convergence rates of Example 4.5  }
\begin{center}
\vspace{-0.1cm}
\begin{tabular}{|c|c|c|c|c|c|c|c|}
\hline
N & \# step  & $\|e_{\vect u_1}\|_{2}$
& rate & $\|e_{\vect u_2}\|_{2}$ & rate &$\|e_{\vect u_3}\|_{2}$ & rate\\
\hline
64   &  13  &1.43e-3  & -        & 1.35e-3  & -        & 2.12e-4 &      -  \\
128   &  13  &3.44e-4  &  2.06 & 2.90e-4  & 2.22  & 3.17e-4 & 2.74 \\
256   &  13  &1.03e-4  &  1.74 & 7.66e-5  & 1.92  & 8.41e-5 & 1.91 \\
512 &  14  &2.91e-5  &  1.82 & 1.99e-5  & 1.94	& 2.15e-5 & 1.97\\	
\hline
\end{tabular}
\end{center}
\label{tabellipse-max}
\end{table}

\begin{table}[!ht]
\caption{$\ell^2$-error and  its  convergence rates of Example 4.5  }
\begin{center}
\vspace{-0.1cm}
\begin{tabular}{|c|c|c|c|c|c|c|}
\hline
N &  $\|e_{\vect u_1}\|_{2}$
& rate & $\|e_{\vect u_2}\|_{2}$ & rate &$\|e_{\vect u_3}\|_{2}$ & rate\\
\hline
64    &8.60e-5  & -        & 1.04e-4  & -        & 6.87e-5 &      -  \\
128    &1.96e-5  &  2.13 & 1.99e-5  & 2.39  & 1.45e-5 & 2.24 \\
256   &4.66e-6  &  1.95 & 4.53e-6  & 2.14  & 3.73e-6 & 1.96 \\
512 &1.19e-6  &  1.97 & 1.15e-6  & 1.98	& 9.63e-7 & 1.95\\	
\hline
\end{tabular}
\end{center}
\label{tabellipse-l2}
\end{table}

\begin{figure}[!ht]
  \centering
    \subfigure[$u_h^{(1)}$]{
    \includegraphics[width=1.5in]{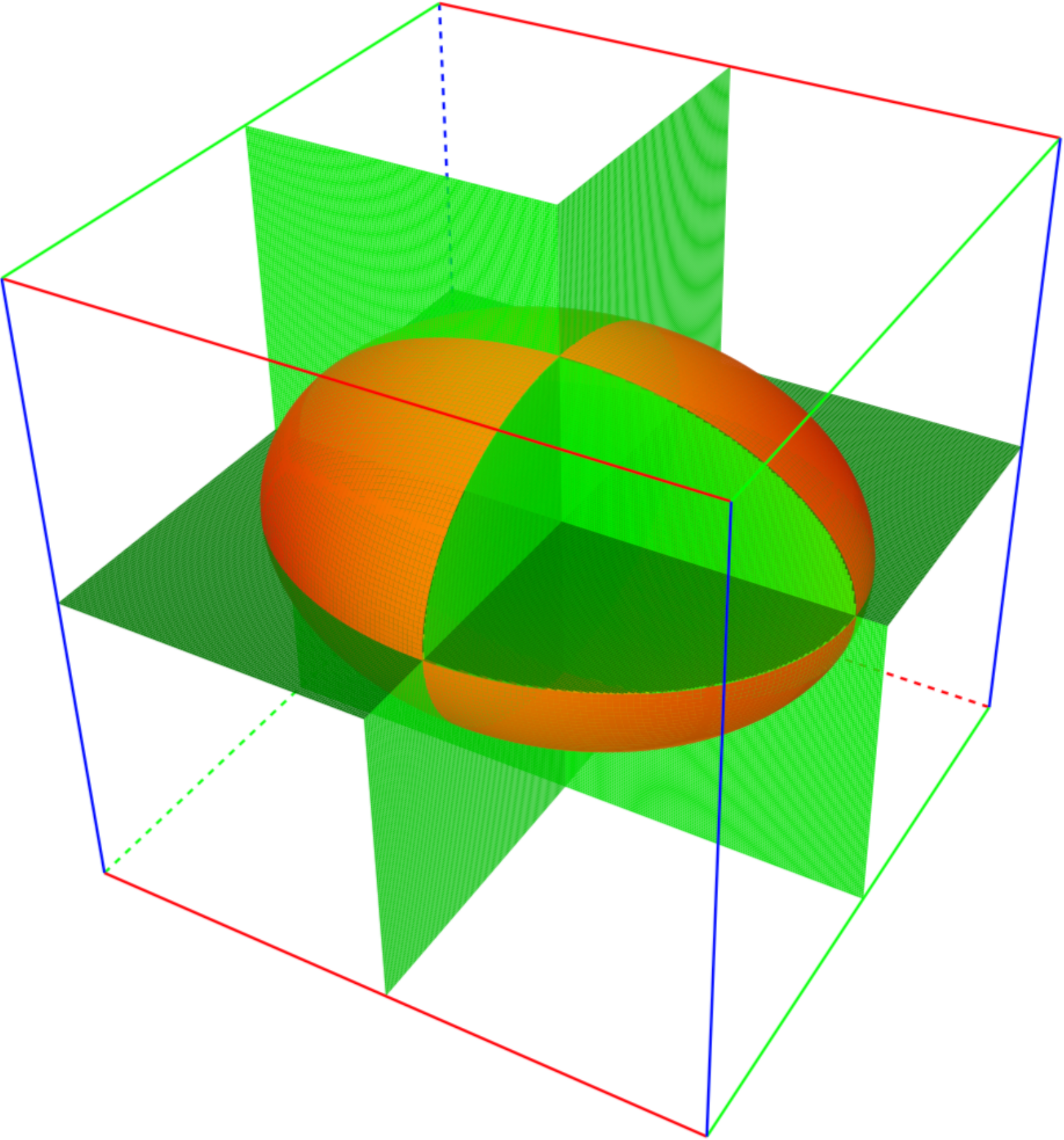}
  }
  \subfigure[$u_h^{(2)}$]{
    \includegraphics[width=1.5in]{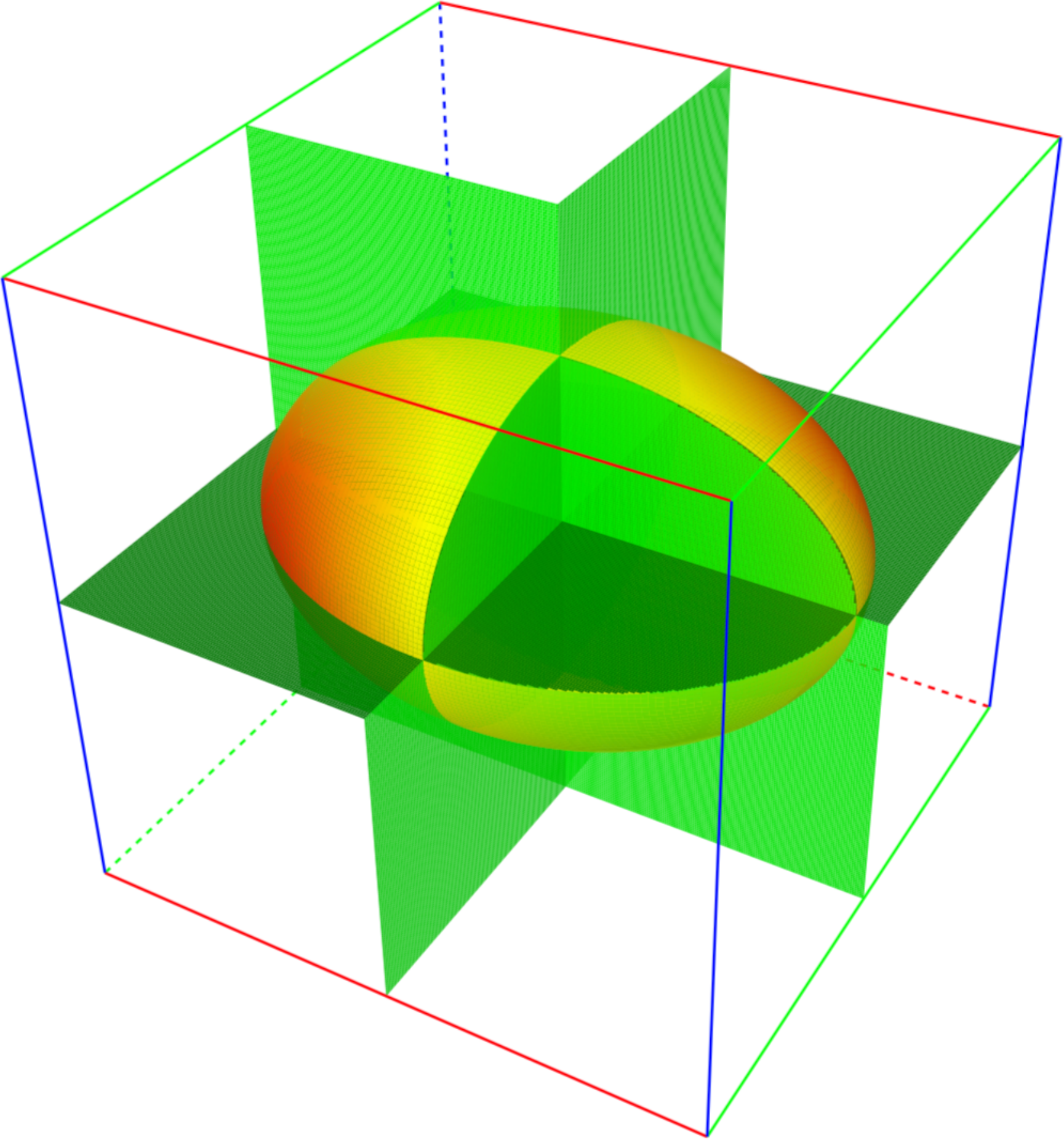}
  }
    \subfigure[$u^{(3)}_h$]{
    \includegraphics[width=1.5in]{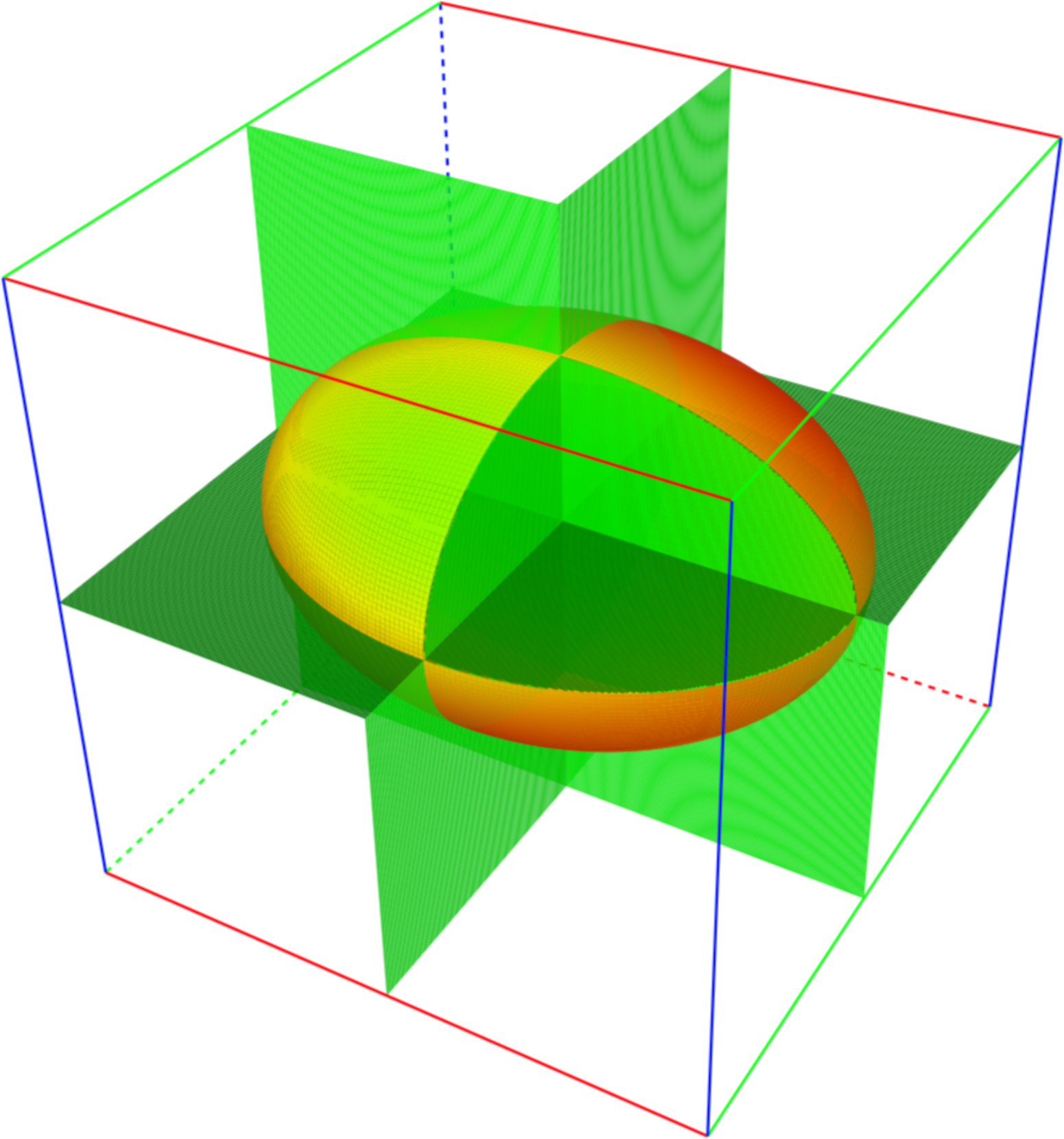}
  }
  \caption{Numerical solution $\vect u_h$ in a ellipsoid for Navier problem.}
    \label{figellipse}
  \end{figure}

{\bf Example 5.6.}
This example solves the problem \eqref{Nm1}-\eqref{Nm2} on a torus which is determined by 
\begin{equation*}
\Omega=\big\{(x,y,z)\in R^3: (c-\sqrt{x^2+y^2})^2+z^2<a^2\big\}
\end{equation*}
with $a=0.35$ and $c=0.7$. The exact solution is chosen by
\begin{equation*}
\begin{cases}
u_1=\exp(1-y^2-z^2)\sin(3x^2+2),\\[1pt]
u_2=\sin(1-x^2-y^2-z^2), \\[1pt]
u_3=3x^2+2y^2+z^2.
\end{cases}
\end{equation*}
The material data are chosen as $E=2000$, and $\nu=2/5$. 
The numerical results in Tables \ref{tabtorus-max} and \ref{tabtorus-l2} verify that the GMRES iteration number is essentially independent of the mesh size and the proposed method also yields second-order accurate solutions. 
\begin{table}[!ht]
\caption{Maximum-error and  its  convergence rates of Example 4.6  }
\begin{center}
\vspace{-0.1cm}
\begin{tabular}{|c|c|c|c|c|c|c|c|}
\hline
N & \# step  & $\|e_{\vect u_1}\|_{2}$
& rate & $\|e_{\vect u_2}\|_{2}$ & rate &$\|e_{\vect u_3}\|_{2}$ & rate\\
\hline
128   &  18  &9.18e-3  &  - & 8.43e-3  & -  & 3.31e-3 & - \\
256   &  19  &1.90e-3  &  2.27 & 1.93e-3  & 2.13  & 7.20e-4 & 2.20 \\
512   &   19 &4.95e-4  &  1.94 & 4.99e-4  & 1.95	& 1.92e-4 & 1.91\\	
\hline
\end{tabular}
\end{center}
\label{tabtorus-max}
\end{table}

\begin{table}[!ht]
\caption{$\ell^2$-error and  its  convergence rates of Example 4.6  }
\begin{center}
\vspace{-0.1cm}
\begin{tabular}{|c|c|c|c|c|c|c|}
\hline
N &  $\|e_{\vect u_1}\|_{2}$
& rate & $\|e_{\vect u_2}\|_{2}$ & rate &$\|e_{\vect u_3}\|_{2}$ &rate\\
\hline
128    &8.74e-4  &  - & 1.08e-3  & -  & 3.30e-4 & - \\
256    &1.77e-4  &  2.30 & 2.34e-4  & 2.21  & 7.22e-5 & 2.19 \\
512    &4.55e-5  &  1.96 & 5.94e-5  & 1.98	& 1.84e-5 & 1.97\\	
\hline
\end{tabular}
\end{center}
\label{tabtorus-l2}
\end{table}

\begin{figure}[!ht]
  \centering
  \subfigure[$u_h^{(1)}$]{
    \includegraphics[width=1.5in]{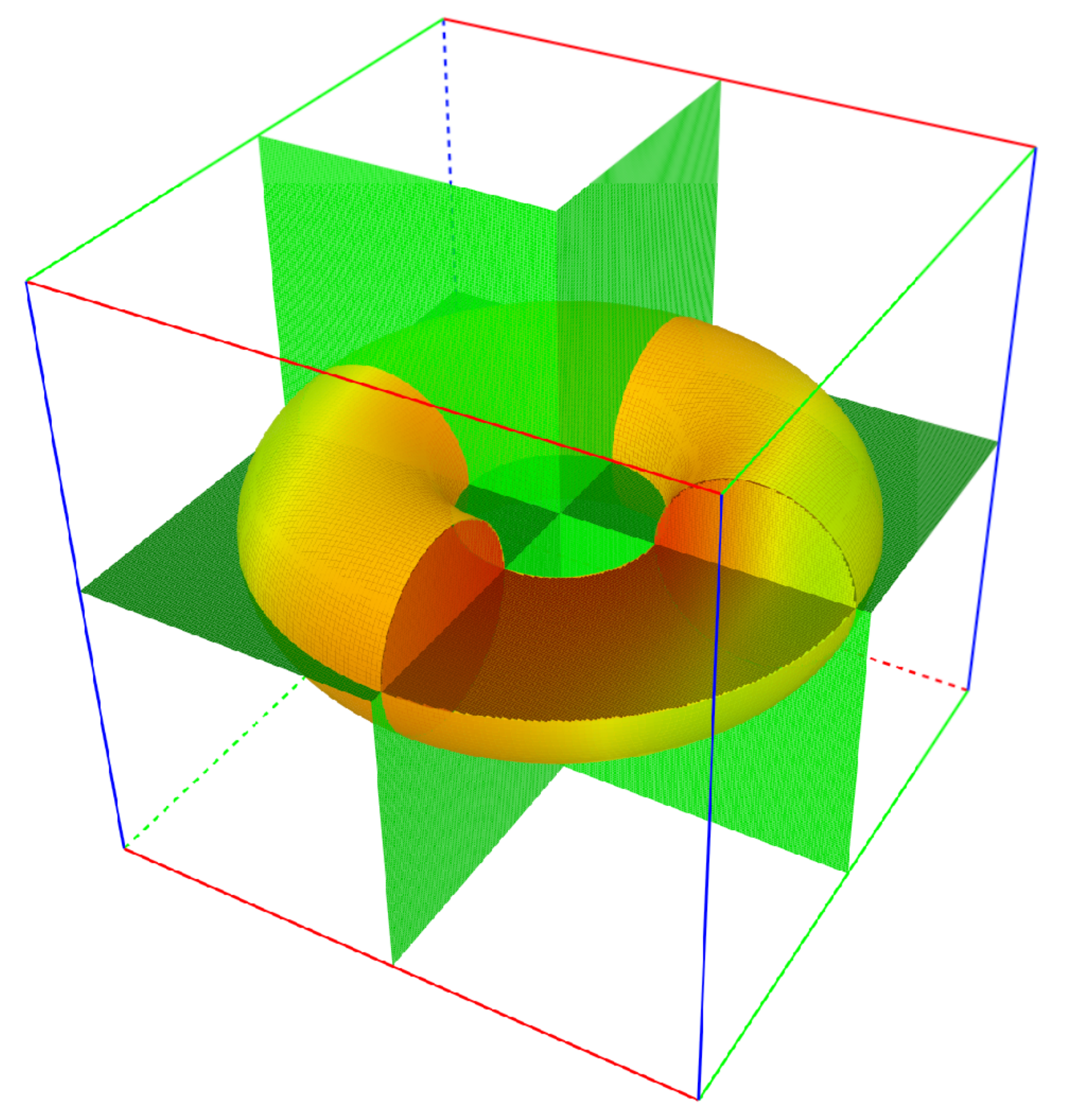}
  }
  \subfigure[$u_h^{(2)}$]{
    \includegraphics[width=1.5in]{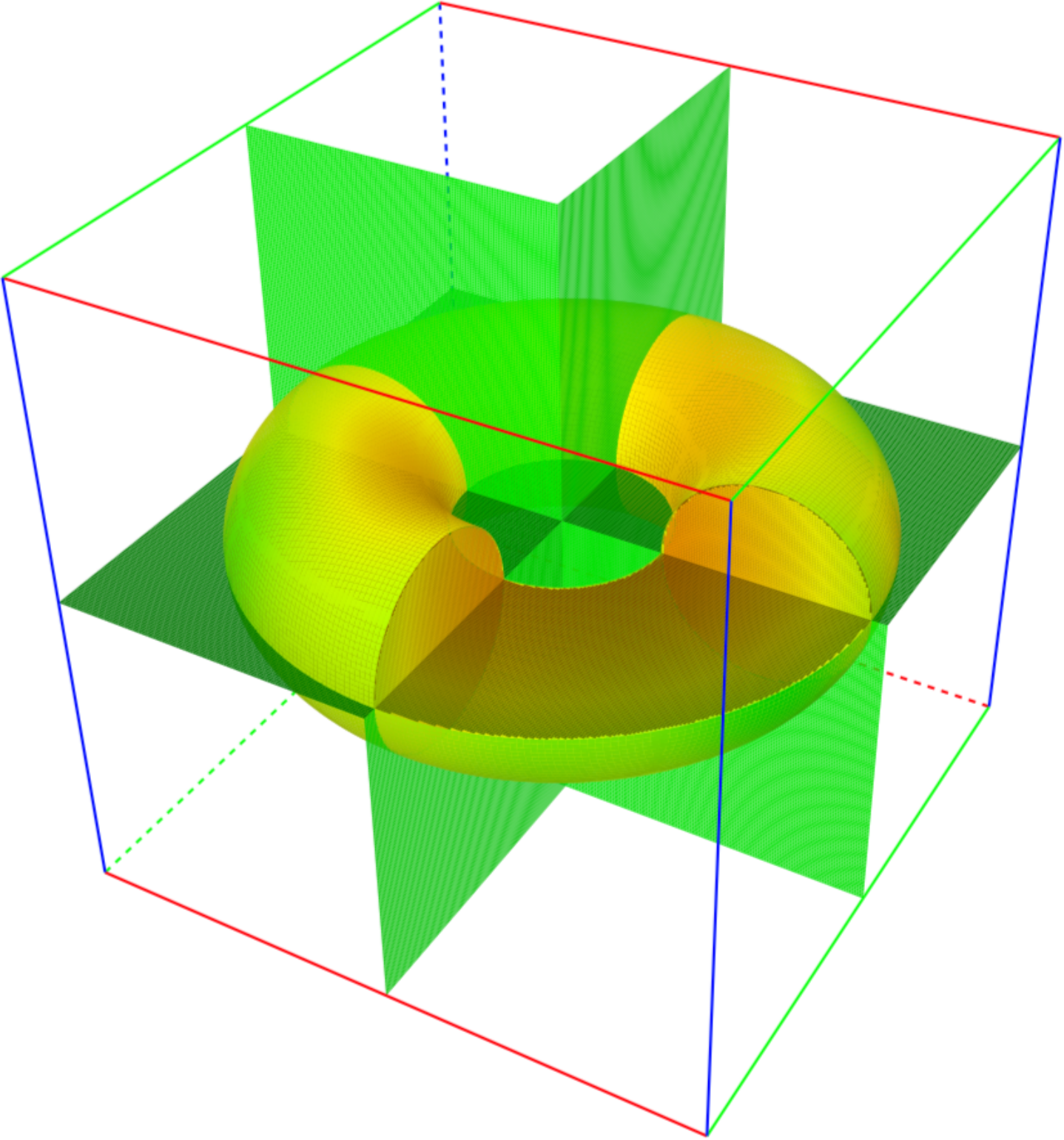}
  }
    \subfigure[$u_h^{(3)}$]{
    \includegraphics[width=1.5in]{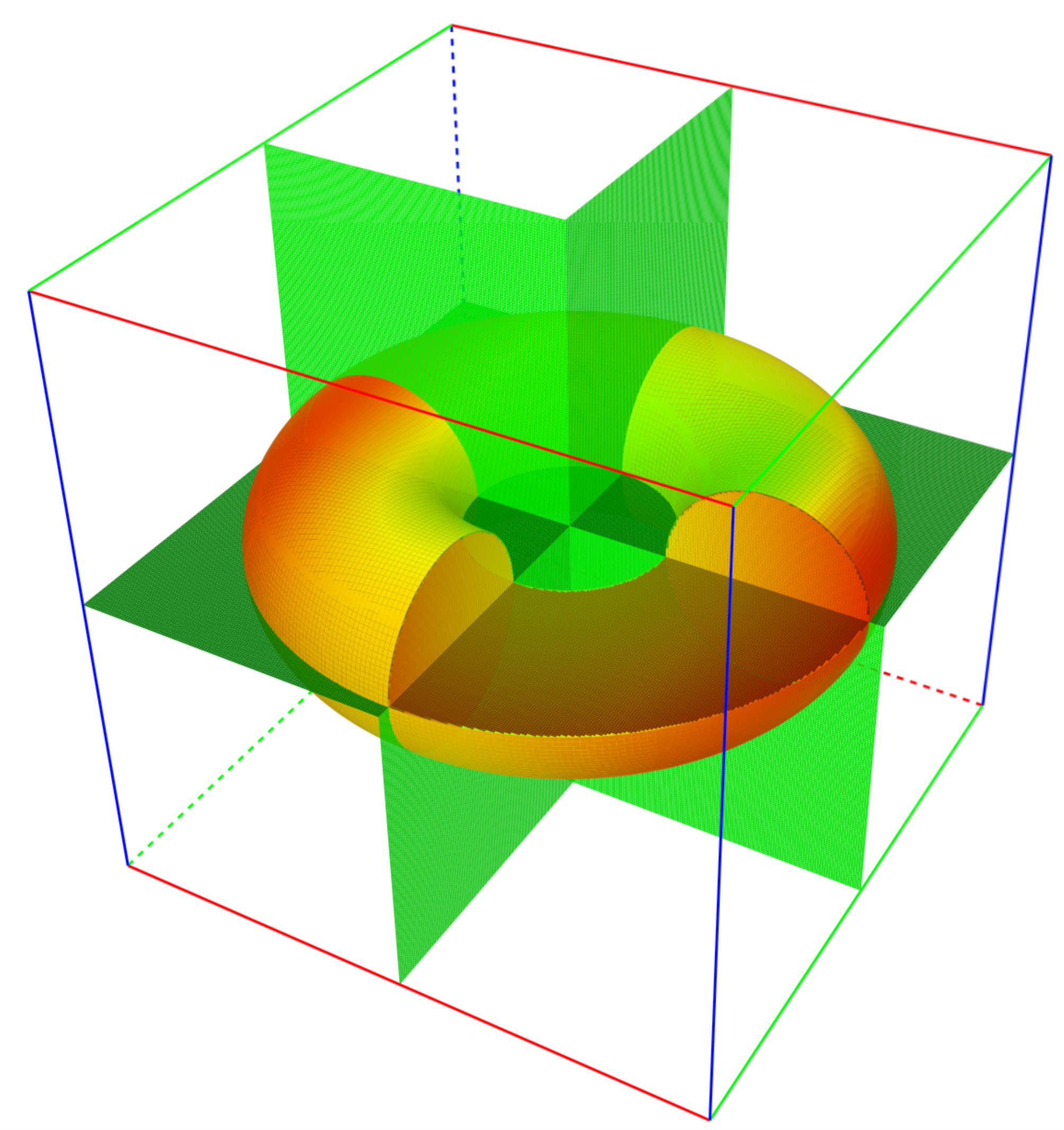}
  }
  \caption{Numerical solution $\vect u_h$ in a torus for Navier problem.}
   \label{figtorus}
  \end{figure}

\section{Conclusions and Discussions}
In this paper, a global second-order KFBI method based on the MAC scheme is proposed to solve the three-dimensional Stokes and Navier boundary value problems on irregular domains. It solves the irregular boundary value problems in the framework of boundary integral equations, but is different from traditional BIM in that the volume or boundary integral is evaluated indirectly. It avoids direct evaluation of nearly singular, singular or hyper-singular boundary integrals, even the requirement of the analytical expressions of Green's functions. 

The reformulated boundary integral equations are Fredholm integral equations of the second kind and can be solved by a Krylov subspace method, such as the GMRES, with the number of iterations being essentially independent of the mesh size. In each GMRES iteration, by posing the Stokes equations in a slightly generalized form that includes a pressure term in the continuity equation, the equivalent simple Stokes and Navier interface problems for integral evaluations can be rewritten into a uniform formulation and  then be  discretized with a modified MAC scheme. Then the discrete system of this scheme is solved efficiently by the CG method together with an FFT-based Poisson solver.  
This approach provides a general algorithmic template for solving two- or multi-fluid problems.  

In addition, the discretization of the surface plays an important role in the KFBI method. This work uses intersection points of the boundary with the grid lines to represent the surface discretization.  The advantage of using intersection points is that it is convenient to find the interpolation stencils, capable to achieve high-order accuracy schemes, and good for problems with moving boundaries.

Nevertheless, the method can be further improved in several aspects. For example,  it suffers deterioration in performance in some cases as the Poisson ratio approaches $1/2$ (i.e., as the material becomes incompressible) for the Navier problem. To solve an almost incompressible elastic material, the technique based on an appropriate decomposition of the Kelvin tensor in \cite{steinbach2007numerical} gives us some ideas. 
This work only describes the details for a second-order version of the KFBI method in three dimensions. In principle, it is natural and straightforward to derive high-order extensions of this method for fluid and solid mechanics. Furthermore, the application of the KFBI method to Stokes-Darcy problems, and Solid-Fluid interaction will be our future work.\\

\section*{Acknowledgments}
Haixia Dong is partially supported by the National Natural Science Foundation of China (Grant No. 12001193), the Scientific Research Fund of Hunan Provincial Education Department (Grant No.20B376), Changsha Municipal Natural Science Foundation (Grant No. kq2014073). Wenjun Ying is partially supported by the Strategic Priority Research Program of Chinese Academy of Sciences (Grant No. XDA25010405), the National Natural Science Foundation of China (Grant No. DMS-11771290) and the Science Challenge Project of China (Grant No. TZ2016002).




\noindent\textbf{References}


\end{document}